\documentclass{amsart}

\usepackage{xcolor}

\usepackage{derivative}
\usepackage{amsthm}
\usepackage{amsmath}
\usepackage{amssymb}
\usepackage{physics}

\usepackage{tikz,pgfplots}
\usepackage{tikz-cd}
\usepackage{float}
\usepackage{stackengine}
\usetikzlibrary{arrows.meta, calc}

\usepackage{pgfplots}
\pgfplotsset{compat=1.18}
\usepackage{booktabs}
\usepackage{multirow}
\usepackage{makecell}
\usepackage{diagbox}

\newtheorem{theorem}{Theorem}[section]

\newtheorem{cor}[theorem]{Corollary}
\newtheorem{prop}[theorem]{Proposition}
\newtheorem{lemma}[theorem]{Lemma}
\newtheorem{defi}[theorem]{Definition}

\theoremstyle{definition}
\newtheorem{rmk}[theorem]{Remark}
\newtheorem{eg}[theorem]{Example}

\numberwithin{equation}{section}

\usepackage{comment}

\title[Classification and Geometry of Loxodromic Curves]{Classification and Geometry of Loxodromic Curves on Surfaces of Revolution in Three-Dimensional Euclidean Space}

\author[Y.C. Huang]{Yen-Chang Huang}
\address{Yen-Chang Huang \\ Department of Applied Mathematics\\National Yang Ming Chiao Tung University, 30010, Taiwan}
\email{ychuang0802@nycu.edu.tw}

\author[Y.H. Chen]{Yan-Hong Chen}
\address{Yan-Hong Chen \\ Department of Applied Mathematics\\National Yang Ming Chiao Tung University, 30010, Taiwan}
\email{cyanhong72@gmail.com}

\thanks{Some of the results presented in this paper were obtained by the second author during undergraduate studies under the supervision of the first author. This work was funded by the National Science and Technology Council (NSTC), Taiwan, under grant numbers 115-2115-M-A49 -001 -MY2 and 112-2115-M-A49 -014 -MY3.
}


\begin{document}

\begin{abstract}
This paper investigates the global geometry of loxodromes
and loxodromic curves on surfaces of revolution in
\(\mathbb{R}^3\).
A loxodrome is a regular curve avoiding the poles and
making a constant directed angle with the meridians,
and it is maximal if it cannot be extended as a
loxodrome on the pole-free surface.
To allow curves to reach or pass through poles, 
we introduce \textit{loxodromic curves}, whose pole-free portions are loxodromes with a common azimuth modulo \(\pi\). A loxodromic curve is maximal if it cannot be extended to a larger parameter interval while remaining a loxodromic curve on the entire surface.

We classify surfaces of revolution into four types:
cylinder-like, torus-like, paraboloid-like, and
sphere-like, and determine the global behavior of
maximal loxodromic curves on each type.
On cylinder-like surfaces, every maximal loxodromic
curve other than a parallel is non-self-intersecting.
On torus-like surfaces, nonparallel maximal loxodromic
curves are closed or dense according to the rationality
of their rotation numbers.
On paraboloid-like and sphere-like surfaces,
only maximal meridians can pass smoothly through poles,
while nonmeridional, nonparallel maximal loxodromic
curves have poles as endpoints.

We derive the differential equations for
arc-length-parametrized loxodromes and introduce
coordinates in which they are represented by straight
lines. We also solve the direct and inverse problems
for loxodromic segments on all four types of surfaces
and characterize maximal loxodromic curves that are
geodesics.

Furthermore, we investigate diffeomorphisms preserving
maximal loxodromic curves and prove that their inverses
have the same preservation property and that the four
surface types are invariant under such diffeomorphisms.
We show that the preservation of maximal meridians and
parallels depends on the surface type: both families
are preserved on paraboloid-like and sphere-like
surfaces, whereas this need not hold on cylinder-like
or torus-like surfaces.
Finally, we characterize azimuth-consistent
diffeomorphisms and establish criteria for
conformality and isometry, including the conformality
of azimuth-consistent diffeomorphisms between
surfaces with poles.

\end{abstract}

\maketitle

\section{Introduction}
A loxodrome, or a rhumb line, on a surface of revolution $S$ is a smooth curve that intersects the meridians at a constant angle (called the azimuth). This concept originated during the Age of Exploration. In 1537, Pedro Nunes distinguished between loxodromes and great-circle routes, and Gerardus Mercator subsequently introduced the Mercator projection in 1569 \cite{alexander2004loxodromes}.

In 1905, Noble provided a parametrization of loxodromes on surfaces of revolution \cite{Noble1905}. In 2009, Kos et al. derived the differential equations for loxodromes based on the first fundamental form \cite{kos2009differential}. In geodesy, where the Earth is commonly modeled by an ellipsoid of revolution, Deakin developed numerical methods for solving the direct and inverse problems for loxodromes on the ellipsoid \cite{deakin2010loxodrome}.

Extensions beyond surfaces of revolution have also been studied. For instance, Babaarslan and Yayli \cite{babaarslan2015differential} extended Noble's classical results \cite{Noble1905} from rotational surfaces to helicoidal surfaces in Euclidean three-space. Caddeo, Onnis, and Piu \cite{caddeo2020loxodromes} considered loxodromes on invariant surfaces in three-dimensional Riemannian manifolds and applied their results to invariant surfaces in $\mathbb{H}^{2}\times\mathbb{R}$, the product of the hyperbolic plane and the real line, and in the three-dimensional Heisenberg group.

In \cite{caddeo2020loxodromes}, the authors derived a local integral representation for loxodromes on invariant surfaces. In the case of Euclidean surfaces of revolution, after adjusting their angle convention and interchanging the profile and orbit coordinates, equation~(2.15) in Remark~2.5 of \cite{caddeo2020loxodromes} becomes
\begin{align}\label{eq123}
   \frac{du}{dv}=\frac{\tan\theta}{r(v)},
\end{align}
which coincides with equation~\eqref{eq:dudv} and is equivalent to the arc-length system in Theorem~\ref{thm:ode_loxo}, with the notation $u,v,\theta,r$ introduced in the next section. They further proved that if a loxodrome making a constant angle $\theta_0\neq0,\pi/2$ with the orbits is also a geodesic, then the invariant surface is flat. In the Euclidean setting, we further investigate the geodesic properties of maximal loxodromic curves on surfaces of revolution. We prove that every maximal loxodromic curve on a circular cylinder is a geodesic (Theorem~\ref{thm:loxodrome_on_a_cylinder}). For surfaces without poles, if the radial function $r$ is not constant on any open interval, a maximal loxodromic curve is a geodesic if and only if it is a maximal meridian or a maximal parallel at a stationary point of $r$ (Theorem~\ref{thm:loxodrome_on_a_no-cylinder}). Moreover, we establish that every maximal meridian is a geodesic on any surface of revolution, including paraboloid-like and sphere-like surfaces where meridians may pass through poles (Proposition~\ref{prop:maximal_meridian_geodesic}).

More recent results on loxodromes are closely related to the present study. Blackwood, Dukehart, and Javaheri \cite{BlackwoodDukehartJavaheri2017} considered a surface of revolution generated by a simple closed profile curve $C(t)=(x(t),y(t))$, $y(t)>0$, parametrized by arc length, and proved that a loxodrome of angle $\psi_0$ is closed if and only if
$$
   \frac{\tan\psi_0}{2\pi}\int_0^T\frac{dt}{y(t)}
   \in\mathbb{Q};
$$
otherwise its trajectory is dense. They also studied the spiraling behavior near a pole. In our notation, their criterion is equivalent to the rationality of $\Omega/2\pi$; see Propositions~\ref{prop:closed_curve} and \ref{prop:dense_loxo_on_torus} below. We place this closed--dense alternative within the classification into cylinder-like, torus-like, paraboloid-like, and sphere-like surfaces developed in Section~\ref{sec:classification}. We also prove that only maximal meridians can pass through a pole; see Corollaries~\ref{cor:topo_of_loxodrome_on_PS} and \ref{cor:topo_of_loxodrome_on_SS}.

Further results on loxodromes have been obtained in several directions. For instance, Kahraman Aksoyak, Bektaş Demirci, and Babaarslan \cite{Aksoyak2023} studied arc-length parametrizations, curvature, and torsion of loxodromes on Euclidean rotational surfaces, while Kovalchuk and Mladenov \cite{KovalchukMladenov2022} solved the direct and inverse loxodromic problems on $\lambda$-spheres. Eastwood \cite{Eastwood2025} studied conformally homogeneous curves on the round two-sphere, with loxodromes among the noncircular cases, and obtained a Möbius-invariant differential equation for their curved analogues. L\'{o}pez \cite{Lopez2026} proved that smooth extremals of a generalized Newton minimal-resistance problem on warped-product surfaces are loxodromes. Since a surface of revolution has metric $dv^2+r(v)^2\,du^2$, the result provides a variational interpretation of the same local constant-angle geometry described by Theorem~\ref{thm:ode_loxo} below.

In this paper, we study the global extensions of loxodromes
and the diffeomorphisms preserving the resulting maximal
loxodromic curves. Particular emphasis is placed on the role
of poles: they constrain both the continuation of curves and
the geometry of curve-preserving maps. In particular, we
prove that azimuth-consistency forces conformality on
surfaces with poles, in contrast to the cylinder-like case.
We now give a section-by-section overview of the main results.

Section~\ref{sec:loxodromes-on-surfaces-of-revolution} introduces the geometric setting and derives the ordinary differential equations for the coordinate functions of an arc-length-parametrized loxodrome; see Theorem~\ref{thm:ode_loxo}. For a nonparallel loxodrome, these equations give \eqref{eq123}. They also lead to the linearizing coordinate $w$, in which loxodromes are represented by straight lines; see Corollary~\ref{cor:uw}.

At a pole, the directed azimuth is undefined. We therefore define a loxodromic curve to be a continuous
curve that is smooth and regular on the interior of its
parameter interval and whose pole-free restrictions are
loxodromes. The directed azimuths of these restrictions agree modulo $\pi$, and their common class in $\mathbb R/\pi\mathbb Z$ is called the azimuth of the loxodromic curve. A maximal loxodromic curve is parametrized by arc length and cannot be extended to a larger interval as an arc-length-parametrized loxodromic curve.

In Section~\ref{sec:classification}, we classify surfaces of revolution $S$ into four types according to the number of poles and whether the profile curve is closed:
\begin{enumerate}
\item cylinder-like (no poles and a nonclosed profile curve),
\item torus-like (no poles and a closed profile curve),
\item paraboloid-like (exactly one pole),
\item sphere-like (exactly two poles).
\end{enumerate}
Here, a pole is a point where the profile curve meets the axis of revolution. Since every pole must correspond to an endpoint of the profile parameter interval, and an interval has at most two endpoints, a surface of revolution has at most two poles under our assumptions. Thus, the number of poles is zero, one, or two. Distinguishing between closed and nonclosed profile curves in the pole-free case yields the four types listed above.

This classification allows us to determine the maximal loxodromic curves on each type of surface. Table~\ref{tab:loxodrome_behaviors_fixed} summarizes whether these curves are closed, non-self-intersecting, or dense and whether they reach or pass through poles. On a torus-like surface, the behavior of a nonparallel maximal loxodromic curve is determined by its rotation number $\Omega/(2\pi)$: a rational value gives a closed curve, whereas an irrational value gives a dense curve. On a surface with poles, only maximal meridians pass through the poles; other nonparallel maximal loxodromic curves may reach poles only as endpoints. The four cases are treated separately in Subsections~\ref{subsec:cylinderlike1}--\ref{subsec:spherelike1}, with corresponding examples provided in each subsection.

  \begin{table}[h!]
   \centering
    \small
    \setlength{\tabcolsep}{4pt}
    \renewcommand{\arraystretch}{1.35}
    \makebox[\textwidth][c]{
        \begin{tabular}{|l|c|p{2.7cm}|p{3.2cm}|l|}
            \hline
            \textbf{Surface Type}
& \makecell[t]{\textbf{Maximal parallel}\\
    $\theta=\pi/2\pmod\pi$}
& \makecell[t]{\textbf{Maximal}\\
    \textbf{meridian}\\
    $\theta=0\pmod\pi$}
& \textbf{Other Azimuths}
& \textbf{Example} \\
            \hline
            Cylinder-like
            & Closed
            & Non-self-intersecting
            & Non-self-intersecting
            & Example \ref{eg:cylinder} \\
            \hline
            Torus-like
            & Closed
            & Closed
            & $\Omega/(2\pi)\in\mathbb{Q}$: closed\newline
              $\Omega/(2\pi)\notin\mathbb{Q}$: non-self-intersecting, dense
            & Example \ref{eg:torus} \\
            \hline
            Paraboloid-like
            & Closed
            & Passes through the pole
            & Non-self-intersecting; pole endpoint
            & Example \ref{eg:hemisphere} \\
            \hline
            Sphere-like
            & Closed
            & Closed and passes through both poles
            & Non-self-intersecting; both poles endpoints
            & Example \ref{eg:sphere} \\
            \hline
        \end{tabular}
    }
    \vspace{5pt}
    \caption{Geometric behavior of maximal loxodromic curves on the four types of surfaces.}
    \label{tab:loxodrome_behaviors_fixed}
\end{table}

In Section~\ref{sec:direct_inverse}, we formulate and solve two classical problems from geodesy on a surface of revolution $S$ of any of the four types: 
\begin{enumerate}
\item \textbf{The direct problem.} Given a nonpole initial point $p_0$ on $S$ and a directed azimuth $\theta$, determine the range of arc lengths $d>0$ for which the corresponding oriented loxodromic segment starting at $p_0$ is defined. For each admissible $d$, determine its terminal point $p_1$.

\item \textbf{The inverse problem.} Given two distinct points $p_0$ and $p_1$ on $S$, determine the possible directed azimuths of the oriented loxodromic segments joining them and their loxodromic distance $d$. Here, the loxodromic distance is defined as the minimum length among all loxodromic segments joining the two points.
\end{enumerate}

For the direct problem, Subsection~\ref{directproblem1} provides the detailed derivation, and Table~\ref{tab:direct_inverse_summary_refs} summarizes the formulas for the four types of surfaces. Each row specifies a surface type, while the columns distinguish meridional, parallel, and all other directed azimuths. The cited equations determine the terminal point $p_1=\Phi(u_1,v_1)\in S$ from the prescribed initial point, directed azimuth, and length, together with the admissible range of lengths. For paraboloid-like and sphere-like surfaces, special formulas are required for meridional segments that may pass through poles, while nonmeridional segments may terminate at a pole but cannot pass through it.

\begin{table}[ht!]
    \centering
    \small
    \setlength{\tabcolsep}{4pt}
    \renewcommand{\arraystretch}{1.35}
    \makebox[\textwidth][c]{
        \begin{tabular}{
        |>{\centering\arraybackslash}m{0.17\textwidth}
        |>{\centering\arraybackslash}m{0.26\textwidth}
        |>{\centering\arraybackslash}m{0.26\textwidth}
        |>{\centering\arraybackslash}m{0.26\textwidth}|}
            \hline
            \textbf{Surface Types}
            & \makecell{\textbf{Meridian}\\
              $\theta=0$ or $\pi$}
            & \makecell{\textbf{Parallel}\\
              $\theta=\pi/2$ or $3\pi/2$}
            & \textbf{Otherwise} \\
            \hline
            Cylinder-like
            & \makecell{
                \eqref{eq:second_point_for_u},
                \eqref{eq:second_point_for_v}
              }
            & \eqref{eq:second_point_for_parl}
            & \makecell{
                \eqref{eq:second_point_for_u},
                \eqref{eq:second_point_for_v}
              } \\
            \hline
            Torus-like
            & \makecell{
                \eqref{eq:second_point_for_u_torus},
                \eqref{eq:second_point_for_v_torus}
              }
            & \eqref{eq:second_point_for_parl_torus}
            & \makecell{
                \eqref{eq:second_point_for_u_torus},
                \eqref{eq:second_point_for_v_torus}
              } \\
            \hline
            Paraboloid-like
            & \makecell{
                \eqref{eq:second_point_for_v_para}, 
                \eqref{eq:second_point_for_u}, 
                \eqref{eq:second_point_for_v}
              }
            & \eqref{eq:second_point_for_parl}
            & \makecell{
                \eqref{eq:second_point_for_u},
                \eqref{eq:second_point_for_v}
              } \\
            \hline
            Sphere-like
            & \makecell{
                \eqref{eq:second_point_for_v_sphere_t},
                \eqref{eq:second_point_for_v_sphere_b}
              }
            & \eqref{eq:second_point_for_parl}
            & \makecell{
                \eqref{eq:second_point_for_u},
                \eqref{eq:second_point_for_v}
              } \\
            \hline
        \end{tabular}
    }
    \vspace{5pt}
    \caption{Summary of equations for the direct problem across different surface types.}
    \label{tab:direct_inverse_summary_refs}
\end{table}

Similarly, Subsection~\ref{inverseproblem1} treats the inverse problem, and Table~\ref{tab:inverse_inverse_summary_refs} lists the formulas used to determine the possible directed azimuths and the loxodromic distance between the prescribed points. The table considers distinct nonpole points $p_0,p_1$. The cases in the columns headed On the Same Parallel and On the Same Maximal Meridian are not mutually exclusive, since two distinct points may lie on both a common parallel and a common maximal meridian. The \textit{Other Cases} column applies when the points lie on neither. For paraboloid-like and sphere-like surfaces, two nonpole points $p_0$ and $p_1$ on the same maximal meridian lie on the same branch if their longitude difference is congruent to $0$ modulo $2\pi$, and on opposite branches if it is congruent to $\pi$ modulo $2\pi$.
\begin{table}[!ht]
    \centering
    \small
    \setlength{\tabcolsep}{4pt}
    \renewcommand{\arraystretch}{1.35}
    \makebox[\textwidth][c]{
        \begin{tabular}{
        |>{\centering\arraybackslash}m{\dimexpr0.16\textwidth-2\tabcolsep-1.25\arrayrulewidth\relax}
        |>{\centering\arraybackslash}m{\dimexpr0.30\textwidth-2\tabcolsep-1.25\arrayrulewidth\relax}
        |>{\centering\arraybackslash}m{\dimexpr0.34\textwidth-2\tabcolsep-1.25\arrayrulewidth\relax}
        |>{\centering\arraybackslash}m{\dimexpr0.25\textwidth-2\tabcolsep-1.25\arrayrulewidth\relax}|}
            \hline
            \textbf{Surface Type}
            & \parbox{\linewidth}{\centering\vspace{3pt} \textbf{On the Same}\\\textbf{Parallel}\par\vspace{3pt}}
            & \parbox{\linewidth}{\centering\vspace{3pt} \textbf{On the Same}\\\textbf{Maximal Meridian}\par\vspace{3pt}}
            & \textbf{Other Cases} \\
            \hline

            Cylinder-like
                & \parbox{\linewidth}{\centering\vspace{3pt} 
                    $\theta=\pi/2$ or $3\pi/2$\\
                    $d_n$: \eqref{eq:distance_of_loxodrome1}
                \par\vspace{3pt}}
                & \parbox{\linewidth}{\centering\vspace{3pt} 
                    $\theta$: \eqref{eq:tangent_of_the_azimuth}\\
                    $d_n$: \eqref{eq:distance_of_loxodrome2}
                \par\vspace{3pt}}
                & \parbox{\linewidth}{\centering\vspace{3pt} 
                    $\theta$: \eqref{eq:tangent_of_the_azimuth}\\
                    $d_n$: \eqref{eq:distance_of_loxodrome2}
                \par\vspace{3pt}} \\
                \hline

            Torus-like
            & \parbox{\linewidth}{\centering\vspace{3pt} 
                $\theta$: $\pi/2$, $3\pi/2$,
                \eqref{eq:tangent_of_the_azimuth_torus}\\
                $d_n$: \eqref{eq:distance_of_loxodrome1}, $d_{n,m}$:\eqref{eq:distance_of_loxodrome3}
              \par\vspace{3pt}}
            & \parbox{\linewidth}{\centering\vspace{3pt} 
                $\theta$: \eqref{eq:tangent_of_the_azimuth_torus}\\
                $d_{n,m}$: \eqref{eq:distance_of_loxodrome3}
              \par\vspace{3pt}}
            & \parbox{\linewidth}{\centering\vspace{3pt} 
                $\theta$: \eqref{eq:tangent_of_the_azimuth_torus}\\
                $d_{n,m}$: \eqref{eq:distance_of_loxodrome3}
              \par\vspace{3pt}} \\
            \hline

            Paraboloid-like
            & \parbox{\linewidth}{\centering\vspace{3pt} 
                $\theta=\pi/2$ or $3\pi/2$\\
                $d_n$: \eqref{eq:distance_of_loxodrome1}
            \par\vspace{3pt}}
            & \begingroup\renewcommand{\arraystretch}{1.1}
\begin{tabular}{@{}>{\centering\arraybackslash}p{\linewidth}@{}}
    \textbf{Same branch}\\
    $\theta$: \eqref{eq:tangent_of_the_azimuth}\\
    $d_n$: \eqref{eq:distance_of_loxodrome2},$d_0$: \eqref{eq:distance_of_same_branch_meridian}\\
    \hline
    \textbf{Opposite branches}\\
    $v_0=v_1$:\\
    $\theta$: $\pi/2$, $3\pi/2$, $0$, $\pi$\\
    $d_n$: \eqref{eq:distance_of_loxodrome1},
    $d_q$: \eqref{eq:distance_of_meridian}\\
    $v_0\neq v_1$:\\
    $\theta$: \eqref{eq:tangent_of_the_azimuth}, $0$ or $\pi$\\
    $d_n$: \eqref{eq:distance_of_loxodrome2},
    $d_q$: \eqref{eq:distance_of_meridian}\\
\end{tabular}\endgroup
            & \parbox{\linewidth}{\centering\vspace{3pt} 
                $\theta$: \eqref{eq:tangent_of_the_azimuth}\\
                $d_n$: \eqref{eq:distance_of_loxodrome2}
              \par\vspace{3pt}} \\
            \hline

            Sphere-like
& \parbox{\linewidth}{\centering\vspace{3pt}
    $\theta$: $\pi/2$, $3\pi/2$\\
    $d_n$: \eqref{eq:distance_of_loxodrome1}
  \par\vspace{3pt}}
& \begingroup\renewcommand{\arraystretch}{1.1}
  \begin{tabular}{@{}>{\centering\arraybackslash}p{\linewidth}@{}}
    \textbf{Same branch}\\
    $\theta$: \eqref{eq:tangent_of_the_azimuth}, $0$, $\pi$\\
    $d_n$: \eqref{eq:distance_of_loxodrome2},
    $d_1, d_2$:\eqref{eq:distance_of_same_longitude_meridians}\\
    \hline
    \textbf{Opposite branches}\\
    $v_0=v_1$:\\
    $\theta$: $\pi/2$, $3\pi/2$, $0$, $\pi$\\
    $d_n$: \eqref{eq:distance_of_loxodrome1},
    $d_a, d_b$: \eqref{eq:distance_of_opposite_longitude_meridians}\\
    $v_0\neq v_1$:\\
    $\theta$: \eqref{eq:tangent_of_the_azimuth}, $0$, $\pi$\\
    $d_n$: \eqref{eq:distance_of_loxodrome2},
    $d_a, d_b$: \eqref{eq:distance_of_opposite_longitude_meridians}
  \end{tabular}\endgroup
& \parbox{\linewidth}{\centering\vspace{3pt}
    $\theta$: \eqref{eq:tangent_of_the_azimuth}\\
    $d_n$: \eqref{eq:distance_of_loxodrome2}
  \par\vspace{3pt}} \\
\hline
        \end{tabular}
    }
    \vspace{5pt}
    \caption{Summary of formulas for the inverse problem between distinct nonpole points.}
    \label{tab:inverse_inverse_summary_refs}
\end{table}

In Section~\ref{sec:geoloxo}, we investigate when maximal loxodromic curves are geodesics on surfaces of revolution. We prove that every maximal loxodromic curve on a circular cylinder is a geodesic (Theorem~\ref{thm:loxodrome_on_a_cylinder}). For surfaces without poles (i.e., cylinder- and torus-like) whose radial function is not constant on any open interval, we show that a maximal loxodromic curve is a geodesic if and only if it is either a maximal meridian or a maximal parallel at a stationary point of $r$ 
(Theorem~\ref{thm:loxodrome_on_a_no-cylinder}). Furthermore, we establish that every maximal meridian is a geodesic on any surface of revolution, including those with poles (Proposition~\ref{prop:maximal_meridian_geodesic}).

Finally, Section~\ref{sec:preserving-maximal-loxodromic-curves} investigates diffeomorphisms $f:S\to\hat S$ between surfaces of revolution that preserve maximal loxodromic curves. Note that this preservation does not require the azimuths to be preserved; the image of a maximal loxodromic curve may have a different azimuth. We prove that the inverse $f^{-1}$ also preserves maximal loxodromic curves and call such a map $f$ a \textit{loxodromic equivalence}. We further show that loxodromic equivalences preserve
both the number of poles and the surface type. We then investigate whether loxodromic equivalences between surfaces of the same type necessarily preserve maximal parallels and maximal meridians. Table~\ref{tab:preservation_properties} summarizes our results.

\begin{table}[ht!]
    \centering
    \begin{tabular}{lcc}
        \toprule
        \textbf{Surface Type} & \makecell{\textbf{Preserves Maximal}\\\textbf{Parallels?}} & \makecell{\textbf{Preserves Maximal}\\\textbf{Meridians?}} \\
        \midrule
        Cylinder-like   & Yes & Not necessarily \\
        Paraboloid-like & Yes & Yes \\
        Sphere-like     & Yes & Yes \\
        Torus-like      & Not necessarily & Not necessarily \\
        \bottomrule
    \end{tabular}
    \vspace{5pt}
    \caption{Preservation properties of loxodromic equivalences on the four types of surfaces.}
    \label{tab:preservation_properties}
\end{table}

We further introduce azimuth-consistent diffeomorphisms, which map all loxodromes with a given directed azimuth on the source surface to loxodromes with a common directed azimuth on the target surface, independently of their initial points. The two azimuths need not coincide. We investigate the differential and geometric properties of azimuth-consistent diffeomorphisms on cylinder-like, paraboloid-like, and sphere-like surfaces, including their conformality and isometry conditions, and study the preservation properties of loxodromic equivalences on torus-like surfaces.

For cylinder-like surfaces, we characterize these diffeomorphisms by the form of their differentials; see Theorem~\ref{thm:LP_and_OGmatrax} and Remark~\ref{rmk:LP_and_Nmatrax}. For a maximal loxodromic curve and its image under an azimuth-consistent diffeomorphism, we also derive an identity relating their respective azimuths $\theta$ and $\hat\theta$ in Remark~\ref{rmk:azimuth_map}. We then establish necessary and sufficient conditions for an azimuth-consistent diffeomorphism to preserve meridians, to be conformal, or to be an isometry; see Corollary~\ref{cor:azimuth_meridian_preservation}, Proposition~\ref{prop:LP_and_Cfm}, and Corollary~\ref{cor:azimuth_isometry}. We also show by Example~\ref{eg:meridian_parallel_not_azimuth_consistent} that a diffeomorphism preserving both meridians and parallels need not be azimuth-consistent and may even fail to preserve maximal loxodromic curves.

For paraboloid-like and sphere-like surfaces, every loxodromic equivalence maps maximal meridians to maximal meridians and maximal parallels to maximal parallels; see Corollary~\ref{cor:map_to_parallel_and_meridian_for_Paraboloid} and Proposition~\ref{prop:map_to_parallel_and_meridian_for_Sphere}. In both cases, we characterize azimuth-consistent diffeomorphisms by the form of their differentials and show that every such diffeomorphism is conformal.

For torus-like surfaces, a loxodromic equivalence need not preserve either maximal meridians or maximal parallels. Using the global linearizing parametrization, we construct such diffeomorphisms from maps of $\mathbb R^2$ associated with matrices in $GL(2,\mathbb Z)$; see Propositions~\ref{prop:f_definition_is_diffeomorphism} and~\ref{prop:loxo_diff_between_tori}. Example~\ref{eg:torus_preservation_examples} gives a map that preserves neither family and a map that preserves both.

\section{Loxodromes on Surfaces of Revolution}\label{sec:loxodromes-on-surfaces-of-revolution}
In this section, we establish the setting and introduce the necessary terminology. Let $S\subseteq \mathbb{R}^3$ be a surface of revolution about the $z$-axis generated by a smooth profile curve
\begin{align}
\alpha(v) = (r(v), 0, z(v)), \label{eq:profile_curve}
\end{align}
which is parametrized by arc length on a nondegenerate interval $I\subseteq\mathbb R$, with $r(v)\geq0$. Let $a=\inf I$ and $b=\sup I$, where $-\infty\leq a<b\leq\infty$. A finite endpoint may or may not belong to $I$, whereas an infinite endpoint never belongs to $I$. A parametrization of $S$ is given by
\begin{align}
\Phi(u, v) = \left(r(v) \cos (u), r(v) \sin (u), z(v)\right), \label{eq:parametrization}
\end{align}
where $u \in \mathbb{R}$ and $v \in I$. A point $q \in S$ is called a \textbf{pole} if $q$ lies on the $z$-axis; namely, $r(v_q) = 0$ for some $v_q \in I$ and $q = \Phi(u, v_q)$ for any $u \in \mathbb{R}$. Throughout the paper, we assume that $S$ is a regular surface, whether or not it has poles. This implies that the tangent vector of the profile curve $\alpha$ at every pole is orthogonal to the $z$-axis; that is, $z'(v_q) = 0$ and $|r'(v_q)|=1$ whenever $r(v_q) = 0$. Note that a surface need not have a pole. For instance, a cylinder with rotation axis along the $z$-axis has no pole since the endpoints of the profile curve (a line segment parallel to the $z$-axis) do not intersect the $z$-axis.

At a finite endpoint \(v_0\in I\), there are two cases: if \(r(v_0)>0\), the parallel
$\{ \Phi(\cdot, v_0)\}$ 
forms a circle of radius $r(v_0)$ and is a boundary component of the surface; if 
$r(v_0)=0$, the corresponding point is a pole.

\begin{rmk}
Under the regularity condition $|r'(v_q)| = 1$, a pole $q = \Phi(u,v_q)$ for any $u\in \mathbb{R}$ occurs only when $v_q$ is at an endpoint of $I$. Indeed, if $v_q$ were an interior point, then $r(v)\geq 0$ and $r(v_q)=0$ would imply that $v_q$ is a minimum point of $r$. Hence $r'(v_q)=0$, a contradiction. Consequently, $S$ has at most two poles, corresponding to the endpoints of $I$. 
\end{rmk}


For a surface of revolution, a \textbf{meridian} (or \textbf{fixed-longitude meridian branch}) is a curve of the form \(\Phi(u_0,v)\), where \(u_0\in\mathbb R\) is fixed and \(v\) varies over \(I\); a \textbf{parallel} is a curve of the form \(\Phi(u,v_0)\), where \(v_0\in I\) is fixed and \(u\) varies over \(\mathbb R\). The partial derivatives of $\Phi$ with respect to $u$ and $v$, which are the respective natural tangent vectors, are
\begin{align}
\partial_u \big\rvert_{\Phi(u,v)} &= (-r(v)\sin u, r(v)\cos u, 0),\label{eq:naturalTV_u} \\
\partial_v\big\rvert_{\Phi(u,v)}&= (r'(v)\cos u, r'(v)\sin u, z'(v))\label{eq:naturalTV_v}
\end{align}
for $u \in \mathbb{R}$ and $v \in I$. Away from the poles, a globally defined orthonormal basis for the tangent spaces of \(S\) is given by
\begin{equation*}
    e_u = \frac{1}{r}\partial_u, \quad e_v = \partial_v.
\end{equation*}

\subsection{Loxodromes} Let $P$ denote the set of poles of $S$, and let $J\subseteq\mathbb R$
be a nondegenerate interval. A regular curve
$\lambda:J\to S\setminus P$ is called a \textbf{loxodrome} if its
unit tangent makes a constant directed angle $\theta\in[0,2\pi)$
with the meridional direction $e_v$; this angle is called its
\textbf{directed azimuth}. For a curve $\lambda:J\to S$ and a point
$p\in S$, we say that $\lambda$ \textbf{passes through} $p$ if
$p\in\lambda(\operatorname{int}J)$. If $p$ is the image of an
endpoint belonging to $J$ and $p\notin\lambda(\operatorname{int}J)$,
we say that $\lambda$ \textbf{reaches} $p$ but does not pass through
it. This distinction is particularly important when $p$ is a pole.


 Let $\lambda:(\epsilon,\delta) \to S$ be a loxodrome of azimuth $\theta$ parametrized by arc length $s$. We have
\begin{equation*}
    \lambda'(s) = \sin(\theta)e_u\rvert_{\lambda(s)} + \cos(\theta)e_v\rvert_{\lambda(s)}.
\end{equation*}
This means that $\lambda$ is an integral curve of the vector field
\begin{equation}\label{eq:Vector_field_loxo_arc}
     W_\theta = \sin(\theta)e_u + \cos(\theta)e_v.
\end{equation}
Note that $W_\theta$ is defined only on $S\setminus P$, since $e_u$ is not
well-defined at the poles. In other words, $\lambda(s) = \Phi(u(s),v(s))$ satisfies the following system of differential equations
\begin{align}
    \frac{du}{ds} &= \frac{\sin\theta}{r(v)}, \label{eq:duds} \\
    \frac{dv}{ds} &= \cos\theta. \label{eq:dvds}
\end{align}
Therefore, we have the following theorem.
\begin{theorem}\label{thm:ode_loxo}
Let $p_0=\Phi(u_0,v_0)\in S\setminus P$, and let $s_0\in(\epsilon,\delta)$. Suppose $\lambda:(\epsilon,\delta)\to S\setminus P$ is a smooth curve parametrized by arc length with $\lambda(s)=\Phi(u(s),v(s))$, where $u(s_0)=u_0$ and $v(s_0)=v_0$. Then $\lambda$ is a loxodrome of directed azimuth $\theta\in[0,2\pi)$ if and only if its coordinate functions satisfy \eqref{eq:duds} and \eqref{eq:dvds}.
\end{theorem}

In addition, integrating these equations yields $u(s)=u_0+\int_{s_0}^s \frac{\sin\theta}{r(v(\tau))}\,d\tau$ and $v(s) = (s-s_0)\cdot\cos(\theta)+v_0$ for initial conditions $u(s_0)=u_0$ and $v(s_0)=v_0$.

A loxodrome may be part of a parallel. To distinguish this case, we call a loxodrome \textbf{nonparallel}, or simply \textbf{not parallel}, if its image is not contained in any parallel. For a loxodrome that is not a parallel, \eqref{eq:dvds} gives $dv/ds=\cos\theta\neq0$. Thus, $v$ is strictly monotonic in $s$, which implies that we can use $v$ to describe the image of the curve on $S$. Writing $s=s(v)$, we can then express the longitude as $u(s(v))$. Applying Theorem~\ref{thm:ode_loxo}, we obtain
\begin{equation}\label{eq:dudv}
        \frac{du}{dv} = \frac{\tan(\theta)}{r(v)}.
\end{equation}
Thus, we may use the surface parameter $v$ to parametrize a nonparallel loxodrome, as stated in the following corollary.

\begin{cor}\label{cor:uv}
    The image of a nonparallel loxodrome of directed azimuth $\theta$ passing through $p_0=\Phi(u_0,v_0)$ is given by
$\Phi(u(v),v)$, where
    \begin{equation}\label{eq:u_of_v}
        u(v) = \tan(\theta) \int_{v_0}^{v}\frac{1}{r(\tau)}\, d \tau + u_0.
    \end{equation}
\end{cor}

We next introduce a new coordinate in which all loxodromes are represented by straight lines. Let $I_0\subseteq I$ be an open interval on which $r(v)>0$, and fix $v_0\in I_0$. Consider the function $w(v) = \int_{v_0}^v \frac{1}{r(\tau)}\, d \tau$ defined on $I_0$. Since $dw/dv=1/r(v)>0$, the function $w$ is strictly increasing and has a smooth inverse $v=v(w)$ on $J=w(I_0)$. Thus, for $v\in I_0$, we obtain the following parametrization
\begin{equation} \label{eq:Conformal_parametrization}
    \Psi(u, w) = \Phi(u, v(w)),
\end{equation}
where $u\in\mathbb{R}$ and $w\in J$. If an endpoint of $I_0$
corresponds to a pole, then $w(v)$ diverges as $v$ approaches
that endpoint. Therefore, this parametrization does not
include the pole itself. The natural tangent vector with respect to $w$ is given by
$$
\partial_w \big\rvert_{\Psi(u,w)} = r(v) \partial_v \big\rvert_{\Phi(u,v)}.
$$
Furthermore, the trajectory of any local loxodrome with azimuth $\theta$ can be parametrized by a curve $\lambda(t)$, which is an integral curve of the vector field
\begin{equation}\label{eq:Vector_field_loxo}
    k(p)\sin(\theta)e_u\rvert_{p} + k(p)\cos(\theta)e_v\rvert_{p},
\end{equation}
for $p\in S\setminus P$ and a positive smooth function $k$ on
$S\setminus P$. Choosing $k(p)=r(v)$ for $p=\Phi(u,v)\in S\setminus P$ and using $\partial_w=r(v)e_v$, we obtain
\begin{align*}
\lambda'(t)
&= r(v)\sin(\theta)e_u\rvert_{\lambda(t)} + r(v)\cos(\theta)e_v\rvert_{\lambda(t)}\\
&= \sin(\theta)\partial_u\rvert_{\lambda(t)} + \cos(\theta)\partial_w\rvert_{\lambda(t)}.
\end{align*}
Thus, we have the following corollary. Unlike the representation in Corollary \ref{cor:uv}, this parametrization applies to all loxodromes, including parallels.

\begin{cor}\label{cor:uw}
    A loxodrome with directed azimuth $\theta$ passing through $p_0=\Phi(u_0,v_0)$ can be expressed as $\lambda(t) = \Psi(u(t),w(t))$, where the coordinate functions are given by
    \begin{align}
    u(t) &= \sin(\theta)t + u_0, \label{eq:ut} \\
    w(t) &= \cos(\theta)t + w_0. \label{eq:wt}
    \end{align}
    Furthermore, if a loxodrome is nonparallel, it can be written as $\lambda(w) = \Psi(u(w),w)$, where
    \begin{equation}\label{eq:uw}
    u(w) = \tan(\theta)(w - w_0) + u_0.
    \end{equation}
\end{cor}

Corollary~\ref{cor:uw} shows that any loxodrome on $S$ can be described as a straight line through the point $(u_0, w_0)$ in the coordinate system $(u, w)$ induced by $\Psi$. In view of this property, we call $w$ the \textbf{linearizing coordinate} for loxodromes and $\Psi$ the \textbf{linearizing parametrization} of $S$. In particular, loxodromes contained in parallels correspond to horizontal lines $w=w_0$, while loxodromes contained in meridians correspond to vertical lines $u=u_0$ in these coordinates. In the subsequent sections, we use this parametrization in the proof of Proposition~\ref{prop:dense_loxo_on_torus} and in the proofs in Subsection~\ref{subsection:Torus-Like Surfaces}.

\subsection{Loxodromic Curves}
Recall that an integral curve of a vector field is \textbf{maximal} if its domain cannot be extended to a larger interval (Definition~14.6 in \cite{Tu2011}); in other words, there is no other integral curve that properly contains it. Accordingly, we say that a loxodrome is \textbf{maximal} if it is a maximal integral curve of a vector field $W_\theta$ defined in \eqref{eq:Vector_field_loxo_arc}.

We are interested in tracing curves on a surface of revolution that maintain a constant angle with the meridians, even when they reach or pass through poles. A natural approach is to consider a maximal solution \((u(s),v(s))\) of \eqref{eq:duds} and \eqref{eq:dvds} and compose it with \(\Phi\), as in Theorem~\ref{thm:ode_loxo}. However, these equations are defined only away from the poles, so the resulting maximal loxodrome may approach a pole without including it.

\begin{figure}
\centering
\includegraphics[scale=0.15]{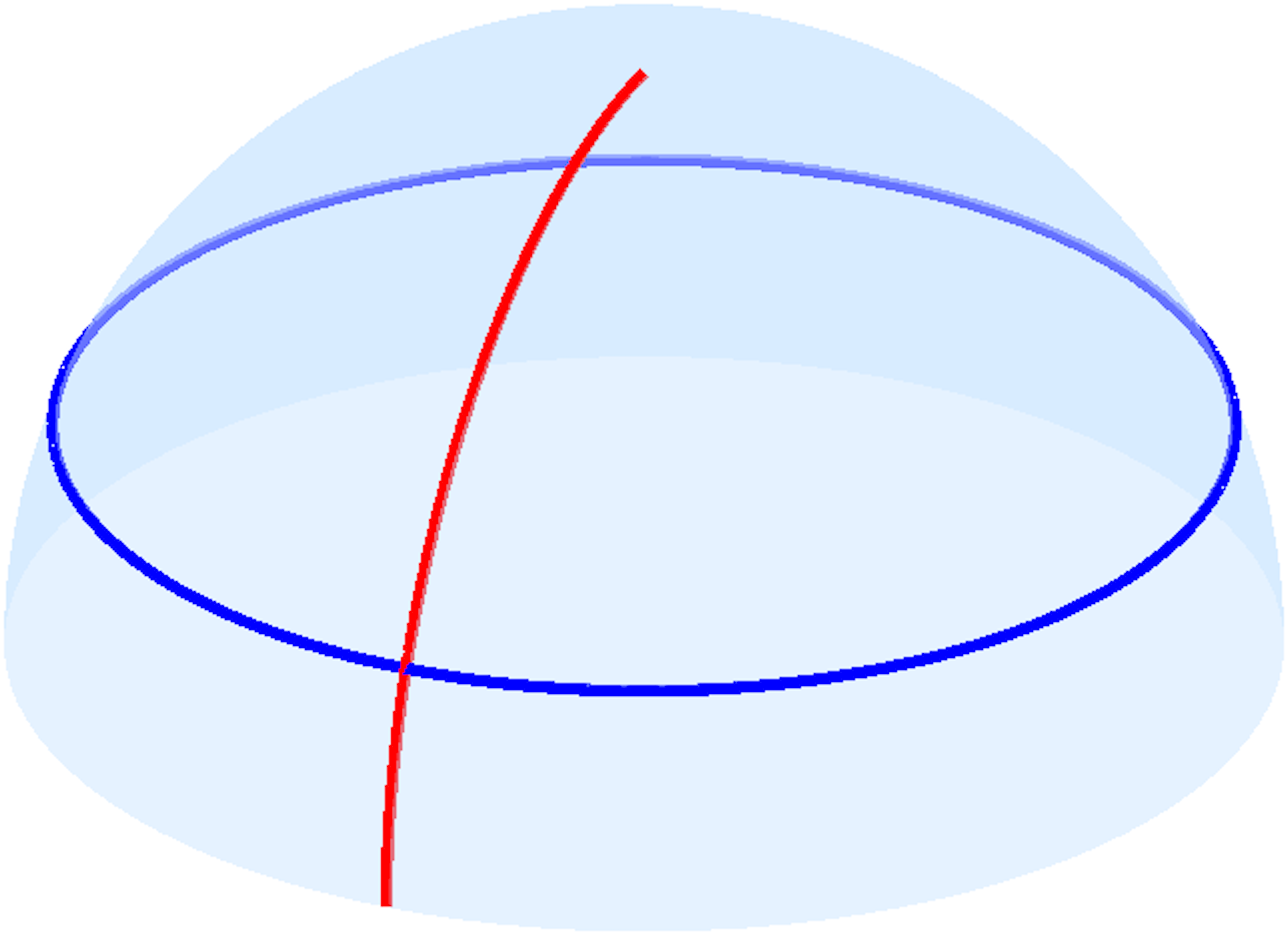}
\caption{\label{fig:parallel_and_meridian} Two maximal loxodromes on the upper unit hemisphere with directed azimuths $\theta=0$ (red) and
$\theta=\pi/2$ (blue). The red curve approaches the north pole but does not include it, while
the blue curve does not pass through any pole.}
\end{figure}

For example, consider the upper unit hemisphere parametrized by
\[
\Phi(u,v)=(\cos v\cos u,\cos v\sin u,\sin v),
\qquad u\in\mathbb R,\quad v\in[0,\pi/2].
\]
The curves
\[
\begin{aligned}
\lambda_1(s)&=(\cos s,0,\sin s),
&&s\in[0,\pi/2),\\
\lambda_2(s)&=\left(
\cos v_0\cos\frac{s}{\cos v_0},
\cos v_0\sin\frac{s}{\cos v_0},
\sin v_0
\right),
&&s\in\mathbb R,
\end{aligned}
\]
where \(v_0\in[0,\pi/2)\) is fixed, are maximal loxodromes of directed azimuths \(\theta=0\) and \(\theta=\pi/2\), respectively, as illustrated in Figure~\ref{fig:parallel_and_meridian}. The red curve \(\lambda_1\) approaches the pole \(p=(0,0,1)\) as \(s\to(\pi/2)^{-}\) but does not include it, whereas the blue curve \(\lambda_2\) is a parallel that avoids the pole. To allow curves to reach or pass through poles, where the directed azimuth is undefined, we introduce the following notion of a loxodromic curve, requiring the constant-angle condition only away from the poles.

\begin{defi}\label{defi:loxodromic_curve}
    Let $J\subseteq\mathbb R$ be a nondegenerate interval. A \textbf{loxodromic curve} is a continuous curve $\lambda : J\to S$ satisfying the following conditions:
    \begin{enumerate}
        \item The curve is smooth and regular on $\operatorname{int}J$.
        \item\label{item:loxodromic-local} For every nonempty open subinterval $J_0\subseteq\operatorname{int}J$ such that $\lambda(J_0)\subseteq S\setminus P$, the restriction $\lambda|_{J_0}$ is a loxodrome.
    \end{enumerate}
\end{defi}

For example, the curve $s\mapsto(\cos s, 0, \sin s)$ for $s\in [0, 2\pi]$ is a loxodromic curve on the unit sphere. 

If $\lambda(\operatorname{int}J)\subseteq S\setminus P$, then, after parametrization by arc length, its coordinate functions satisfy \eqref{eq:duds} and \eqref{eq:dvds} throughout the interior of its parameter domain by Theorem~\ref{thm:ode_loxo}. The following lemma shows that the directed azimuths of these restrictions agree modulo $\pi$, even when the curve passes through a pole.

\begin{lemma}\label{lemma:modulo_azimuths}
Let $\lambda : J\to S$ be a loxodromic curve. Then the directed azimuths of all restrictions $\lambda|_{J_0}$ described in condition~\ref{item:loxodromic-local} of Definition \ref{defi:loxodromic_curve} agree modulo $\pi$.
\end{lemma}

The proof of this lemma will be given case by case in the next section. The common class of these directed azimuths in $\mathbb R/\pi\mathbb{Z}$ is called the \textbf{azimuth} of $\lambda$. 

An arc-length-parametrized loxodromic curve $\lambda:K\to S$ is called \textbf{maximal} if there do not exist an interval $J$ properly containing $K$ and an arc-length-parametrized loxodromic curve $\mu:J\to S$ such that $\mu|_K=\lambda$.




\section{Classification of Surfaces of Revolution and Behavior of Loxodromes}\label{sec:classification}

In this section, we classify surfaces of revolution according to the number of poles and whether the profile curve is closed. Depending on the number of poles and the behavior of the profile curve, we establish the following classification: if $S$ has no poles and its profile curve is nonclosed, it corresponds to the standard case commonly discussed in classical differential geometry, which we call a \textbf{cylinder-like surface}; if $S$ has no poles and its profile curve is a closed curve, we call it a \textbf{torus-like surface}; if $S$ has exactly one pole, we call it a \textbf{paraboloid-like surface}; finally, if $S$ has exactly two poles, we call it a \textbf{sphere-like surface}.

\subsection{Cylinder-Like Surfaces}\label{subsec:cylinderlike1}
A cylinder-like surface $S$ is a surface without poles, generated by a non-self-intersecting, nonclosed profile curve. Note that, by the height of $S$, we mean the value $\sup_{v_1,v_2\in I}|z(v_1)-z(v_2)|$. Since $S$ has no poles, we have $r(v) > 0$ for all $v \in I$.
\begin{proof}[Proof of Lemma~\ref{lemma:modulo_azimuths} for cylinder-like surfaces]
Let $J_1,J_2\subseteq\operatorname{int}J$ be nonempty open subintervals. Since $S$ has no poles, the restriction of the loxodromic curve $\lambda$ to any open subinterval of $\operatorname{int}J$ is a loxodrome.

If $J_1\cap J_2\neq\varnothing$, the restrictions $\lambda|_{J_1}$ and $\lambda|_{J_2}$ coincide on their overlap and therefore have the same directed azimuth.

If $J_1\cap J_2=\varnothing$, choose an open subinterval $J_0\subseteq\operatorname{int}J$ that intersects both $J_1$ and $J_2$. By the preceding argument, each of $\lambda|_{J_1}$ and $\lambda|_{J_2}$ has the same directed azimuth as $\lambda|_{J_0}$. Thus, their directed azimuths agree, and in particular agree modulo $\pi$.
\end{proof}
\begin{rmk}\label{rmk:same_directed_azimuth_cylinder}
    By the proof of Lemma~\ref{lemma:modulo_azimuths}  for cylinder-like surfaces, all restrictions of $\lambda$ to nonempty open subintervals of $\operatorname{int}J$ have the same directed azimuth.
\end{rmk}
\begin{lemma}\label{lemma:maximal_solution_represents_loxodromic_curve} 
    Suppose $S$ is a cylinder-like surface. Let $\lambda:J\to S$ be a maximal arc-length-parametrized loxodromic curve, and let $\theta$ be the common directed azimuth of its restrictions. If $s_0\in\operatorname{int}J$ and $p_0=\lambda(s_0)=\Phi(u_0,v_0)$, then the maximal solution $(u,v)$ of \eqref{eq:duds} and \eqref{eq:dvds} satisfying $u(s_0)=u_0$ and $v(s_0)=v_0$ has domain $\operatorname{int}J$ and satisfies $\lambda(s)=\Phi(u(s),v(s))$ for every $s\in\operatorname{int}J$.
\end{lemma}

\begin{proof}
Since $S$ has no poles and $\operatorname{int}J$ is an interval, we can choose  coordinate functions $\hat u,\hat v$ along $\lambda$ on $\operatorname{int}J$ such that $\lambda(s)=\Phi(\hat u(s),\hat v(s))$, where $\hat u(s_0)=u_0$ and $\hat v(s_0)=v_0$. By Remark~\ref{rmk:same_directed_azimuth_cylinder} and Theorem~\ref{thm:ode_loxo}, $(\hat u,\hat v)$ satisfies \eqref{eq:duds} and \eqref{eq:dvds} on $\operatorname{int}J$.

Let $(u,v)$ be the maximal solution on an open interval $K$ with the same initial conditions, and define $\mu(s)=\Phi(u(s),v(s))$ for each $s\in K$. By the existence and uniqueness theorem for ordinary differential equations, $K$ contains $\operatorname{int}J$ and $(u,v)|_{\operatorname{int}J}=(\hat u,\hat v)$. Therefore, $\mu|_{\operatorname{int}J}=\lambda|_{\operatorname{int}J}$.

Suppose that $K$ properly contains $\operatorname{int}J$. Since $K$ is an open interval, $J\cup K$ properly contains $J$. Define $\nu:J\cup K\to S$ by
$$
\nu(s)=
\begin{cases}
\lambda(s), & \text{if } s\in J,\\
\mu(s), & \text{if } s\in K\setminus J.
\end{cases}
$$
The curves $\lambda$ and $\mu$ agree on $\operatorname{int}J$. By continuity, they also agree at every endpoint of $J$ contained in $K$. Hence $\nu$ is well-defined. Moreover, near every endpoint of $J$ contained in the interior of $J\cup K$, the curve $\nu$ agrees with $\mu$. Therefore, $\nu$ is an arc-length-parametrized loxodromic curve. Since $
\nu|_J=\lambda$  and $J\cup K$ properly contains $J$, this contradicts the maximality of $\lambda$. Therefore, $K=\operatorname{int}J$.
\end{proof}

 Therefore, by Lemma~\ref{lemma:maximal_solution_represents_loxodromic_curve}, Theorem~\ref{thm:ode_loxo} applies throughout the interior of the domain of every maximal loxodromic curve on a cylinder-like surface.

\begin{figure}
\centering
\includegraphics[scale=0.3]{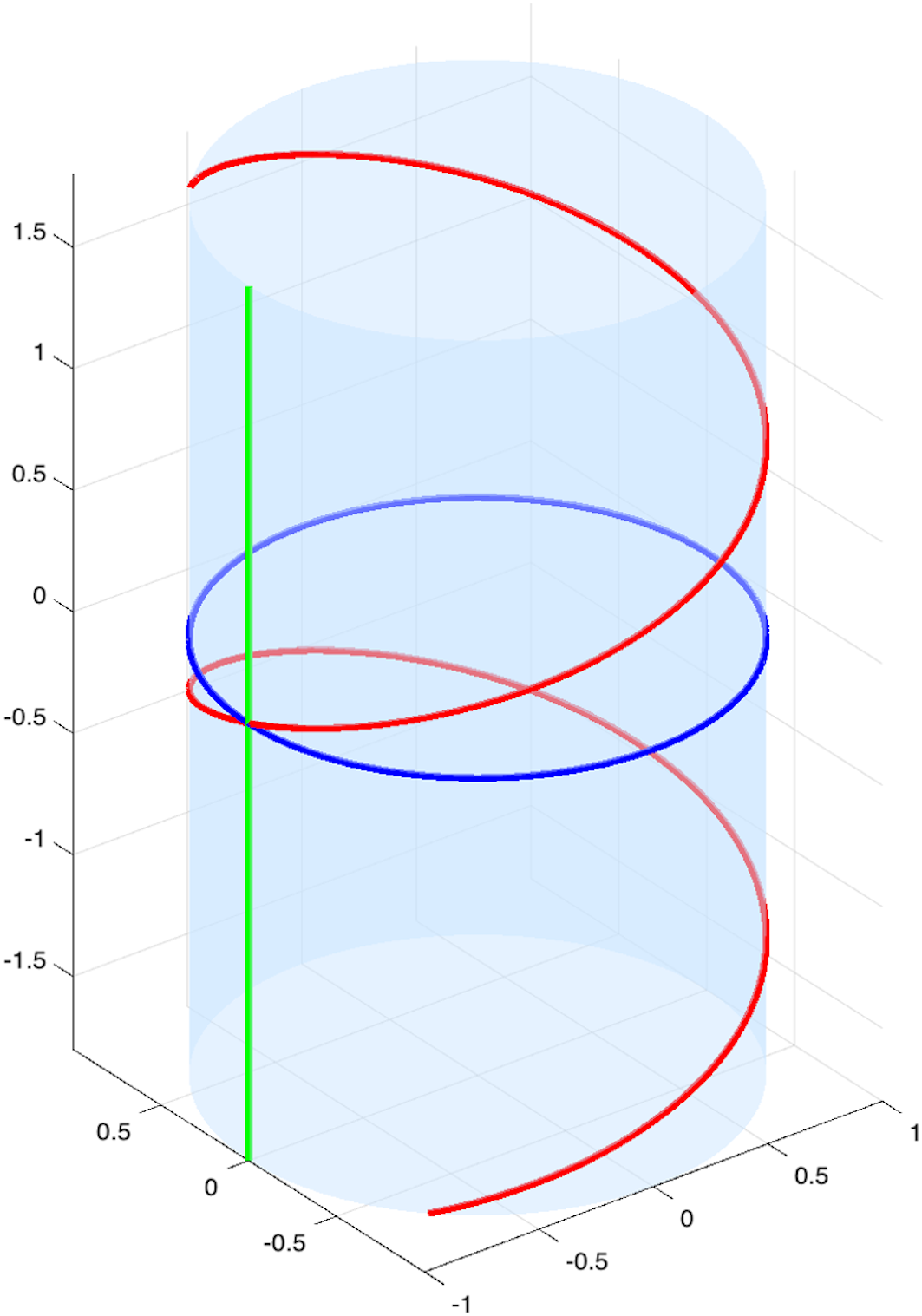}
\caption{\label{fig:cylinder} A loxodromic  helix (red), a meridian (green), and a parallel (blue) on the cylinder.}
\end{figure}

\begin{prop}\label{prop:topo_of_loxodrome_on_CS}
    Let $S$ be a cylinder-like surface. There are only two kinds of maximal loxodromic curves:
    \begin{enumerate}
        \item maximal parallels, and
        \item non-self-intersecting maximal loxodromic curves with finite or infinite length, depending on the arc length of the profile curve.
    \end{enumerate}
\end{prop}
Here, \textbf{maximal parallels} are maximal loxodromic curves with azimuth $\pi/2$.
\begin{proof}
     Let $\lambda:J\to S$ be a maximal loxodromic curve. If the azimuth of $\lambda$ is $\pi/2 \bmod \pi$, then it is a maximal parallel. Otherwise, let $\theta$ be the common directed azimuth of its restrictions. By Lemma~\ref{lemma:maximal_solution_represents_loxodromic_curve}, $\lambda$ is given by $\lambda(s)=\Phi(u(s),v(s))$, where $(u,v)$ is the maximal solution of \eqref{eq:duds} and \eqref{eq:dvds}. By continuity, the formula for $v$ also holds at every included endpoint of $J$. Since $\cos\theta\neq0$, the function $v$ is injective. Let $s_1,s_2$ belong to the parameter interval of $\lambda$, and suppose that $\lambda(s_1)=\lambda(s_2)$. By \eqref{eq:parametrization}, we have $z(v(s_1)) = z(v(s_2))$, along with
\begin{align*}
    r(v(s_1))\cos(u(s_1)) &= r(v(s_2))\cos(u(s_2)), \\
    r(v(s_1))\sin(u(s_1)) &= r(v(s_2))\sin(u(s_2)).
\end{align*}
Squaring and adding the last two equations yields $r(v(s_1)) = r(v(s_2))$. Because the profile curve is non-self-intersecting, we have $v(s_1) = v(s_2)$. By the injectivity of $v$, we conclude that $s_1 = s_2$. Thus, $\lambda$ is injective. By the maximality of $(u,v)$ and the positivity of $r$ on $I$, the image of $v$ on $\operatorname{int}J$ is $\operatorname{int}I$. Since $dv/ds=\cos\theta\neq0$, the arc length of $\lambda$ is $|I|/|\cos\theta|$ by \eqref{eq:dvds}. Consequently, its length is finite if and only if the profile curve has finite length.
\end{proof}
\begin{eg}\label{eg:cylinder}
    Given a profile curve $\alpha(v) = (1,0,v)$ on the interval $I$, a parametrization of $S$ is
    $$
    \Phi(u,v) = (\cos(u),\sin(u),v).
    $$
    In particular, when $|I| < \infty$, $S$ is of finite height; otherwise, $S$ has infinite height. By Lemma~\ref{lemma:maximal_solution_represents_loxodromic_curve}, the maximal loxodromic curve of azimuth $\theta$ passing through $p_0 = \Phi(u_0,v_0)$  is $\lambda(s)=\Phi(u(s),v(s))$, where
    \begin{align*}
        u(s) &=\pm \sin(\theta) \cdot (s-s_0)+u_0,\\
        v(s) &= \pm \cos(\theta) \cdot (s-s_0)+v_0.
    \end{align*}
    Here, $\theta\in[0,\pi)$ represents the azimuth class in $\mathbb R/\pi\mathbb Z$, and the same sign is chosen
in both equations. As illustrated in Figure~\ref{fig:cylinder}, the red curve on this cylinder forms a helix with a constant slope.
\end{eg}

Let $\lambda:J\to S$ be a nonparallel maximal loxodromic curve, and let $\theta$ be the common directed azimuth of its restrictions to nonempty open subintervals of $\operatorname{int}J$. Suppose that $\lambda(s_0)=\Phi(u_0,v_0)$. Since $v(s)=\cos\theta(s-s_0)+v_0$, Lemma~\ref{lemma:maximal_solution_represents_loxodromic_curve} gives
$$
\operatorname{int}J=
\begin{cases}
\left(s_0+\sec\theta(a-v_0),\,s_0+\sec\theta(b-v_0)\right),
& \text{if } \cos\theta>0,\\
\left(s_0+\sec\theta(b-v_0),\,s_0+\sec\theta(a-v_0)\right),
& \text{if } \cos\theta<0.
\end{cases}
$$
Using $v$ as the parameter, we can write $\lambda$ in the form given in Corollary~\ref{cor:uv} and \eqref{eq:u_of_v}, where $v\in\operatorname{int}I$. Using $w$ as the parameter, we can also write $\lambda$ in the form given in Corollary~\ref{cor:uw} and \eqref{eq:uw}, where $w\in w(\operatorname{int}I)$.
\begin{rmk}
Suppose the profile curve $\alpha:I \rightarrow \mathbb{R}^3$ in this case is an embedding. Then $S$ is homeomorphic to $I\times \mathbb S^1$, where
$\mathbb S^1=\{e^{2\pi i x}\in\mathbb{C}:x\in\mathbb{R}\}$.
More precisely, the homeomorphism is given by
$$
(v,e^{i\theta})\longmapsto \Phi(\theta,v),
$$
where $v\in I$ and $\theta\in\mathbb{R}$. This map is well defined because $\Phi$ is $2\pi$-periodic in its first variable. Thus, every cylinder-like surface is homeomorphic to $I\times\mathbb S^1$, and hence its fundamental group is isomorphic to $\mathbb Z$. If the embedding assumption on the profile curve is removed, the resulting revolved set need not be a regular surface. For example, consider the injective profile curve
$$
\alpha(v) = (\sin(2v)+2, 0, \cos v)
$$
for each $v \in (-3\pi/2,\pi/2)$. Its image is homeomorphic to the wedge sum $S^1\vee S^1$. The set obtained by revolving the image of $\alpha$ is not a regular surface and is therefore excluded from our classification. Nevertheless, since $r(v)>0$, this revolved set is homeomorphic to $(S^1\vee S^1)\times S^1$. Consequently, its fundamental group is $(\mathbb Z*\mathbb Z)\times\mathbb Z$, where $*$ denotes
the free product of groups. This fundamental group is nonabelian.
\end{rmk}

\subsection{Torus-Like Surfaces}\label{subsec:toruslike1}
Let $S$ be a torus-like surface generated by a profile curve $\alpha:[a,b] \to \mathbb{R}^3$ parametrized by arc length. By the definition of a torus-like surface, we have $\alpha(a)=\alpha(b)$ and $\alpha'(a)=\alpha'(b)$. Since $S$ has no poles, we have $r(v) > 0$ for all $v \in [a,b]$. In \cite{Pressley2010}, the Gaussian curvature of a surface of revolution is given by $K = \frac{-r''}{r}$. Using the curvature formula, we obtain
\begin{equation}\label{eq:G-B_Thm}
    \int_S K \,dA =\int_0^{2\pi}\int_a^b \bigg(\frac{-r''(v)}{r(v)}\bigg)r(v)\,dv\,du=2\pi (r'(a)-r'(b)).
\end{equation}
For a torus-like surface, the periodicity of the profile curve
implies $r'(a)=r'(b)$. Hence $\int_S K\,dA=0$.
Since $S$ is a compact, connected, orientable surface without
boundary, the Gauss--Bonnet theorem gives $\chi(S)=0$.
By the classification of compact connected orientable surfaces, $S$ is homeomorphic to a torus $\mathbb{T}^2=\mathbb{S}^1\times\mathbb{S}^1$. Consequently, its fundamental group is $\mathbb Z\times\mathbb Z$, which is abelian.

The proof of Lemma~\ref{lemma:modulo_azimuths} for torus-like surfaces is the same as that for cylinder-like surfaces. Thus, all restrictions of a loxodromic curve to nonempty open subintervals of the interior of its domain have the same directed azimuth. To analyze the global behavior of maximal loxodromic curves, we extend $\alpha$ to a periodic function $\tilde{\alpha} = (\tilde{r}, 0, \tilde{z})$ with period $L = b - a$ on the entire real line $\mathbb{R}$; that is, for each $n \in \mathbb{Z}$ and $v \in [a+nL, b+nL]$, we have
\begin{equation}\label{extended_profile_curve}
    \tilde{\alpha}(v) = \alpha(v-nL).
\end{equation}
Using \eqref{extended_profile_curve}, we describe the torus-like surface by defining a new parametrization
\begin{equation*}
    \tilde{\Phi}(u,v) = (\tilde{r}(v)\cos(u), \tilde{r}(v)\sin(u), \tilde{z}(v))
\end{equation*}
for each $u,v \in \mathbb{R}$. This parametrization satisfies $\tilde{\Phi}(u,v) = \Phi(u, v-nL)$ for any $n \in \mathbb{Z}$ and $v \in [a+nL, b+nL]$. 

\begin{lemma}\label{lemma:maximal_solution_torus}
Let $S$ be a torus-like surface, and let $\lambda:J\to S$
be a maximal loxodromic curve parametrized by arc length.
Let $\theta\in[0,2\pi)$ denote the directed azimuth of
$\lambda|_{\operatorname{int}J}$.
Fix $s_0\in\operatorname{int}J$, and write
$p_0=\lambda(s_0)=\Phi(u_0,v_0)$.
Then the maximal solution $(\tilde u,\tilde v)$ of
\begin{align}
    \frac{d\tilde u}{ds}
    =\frac{\sin\theta}{\tilde r(\tilde v)},
    \qquad
    \frac{d\tilde v}{ds}=\cos\theta,
\end{align}
with initial conditions
$\tilde u(s_0)=u_0$ and $\tilde v(s_0)=v_0$,
has domain $\operatorname{int}J$ and satisfies
$\lambda(s)=\tilde\Phi(\tilde u(s),\tilde v(s))$
for every $s\in\operatorname{int}J$.
\end{lemma}

The proof is similar to that of Lemma~\ref{lemma:maximal_solution_represents_loxodromic_curve}. Since $\operatorname{int}J$ is an interval, we can choose smooth real-valued coordinate functions along $\lambda$ with the prescribed initial values, using $\tilde{\Phi}$. 


The maximal loxodromic curves with azimuth $\frac{\pi}{2} \bmod \pi$ are maximal parallels. We now consider nonparallel maximal loxodromic curves.
Let $\lambda:J\to S$ be a nonparallel maximal loxodromic curve, and let $\theta$ be the common directed azimuth of its restrictions to nonempty open subintervals of $\operatorname{int}J$. Suppose that $\lambda(s_0)=\tilde{\Phi}(u_0,v_0)$. Since $\tilde{v}(s)=\cos\theta(s-s_0)+v_0$, Lemma~\ref{lemma:maximal_solution_torus} gives $\operatorname{int}J=(-\infty,\infty)$. Using $\tilde v$ as the parameter, we can express $\lambda$ as $\tilde\Phi(\tilde{u}(\tilde{v}),\tilde{v})$ with
\begin{equation}\label{eq:u_of_v_tilde}
        \tilde u(\tilde v) = \tan(\theta) \int_{v_0}^{\tilde v}\frac{1}{\tilde r(\tau)}\, d \tau + u_0,
    \end{equation}
which is given by (\ref{eq:u_of_v}) in Corollary~\ref{cor:uv}, where $\tilde v\in\mathbb R$. Since $\tilde r$ is periodic with period $L=b-a$, for any $\tilde v \in \mathbb{R}$, we have
\begin{equation}\label{eq:quasi-periodic_of_u}
    \tilde u(\tilde v+L) - \tilde u(\tilde v) = \tan(\theta) \int_{\tilde v}^{\tilde v+L} \frac{1}{\tilde{r}(\tau)} \,d\tau = \tan(\theta) \int_{a}^b \frac{1}{r(\tau)} \,d\tau.
\end{equation}
Define
\begin{equation}\label{eq:periodic_u}
    \Omega = \tan(\theta) \int_{a}^b \frac{1}{r(\tau)} \,d\tau.
\end{equation}
For fixed $\theta$, we have $\tilde u(\tilde v+L)=\tilde u(\tilde v)+\Omega$
for every $\tilde v\in\mathbb{R}$. Thus, $\Omega$ is the difference between the longitudes $\tilde u(\tilde v+L)$ and $\tilde u(\tilde v)$, corresponding to one full circuit of the profile curve.

\begin{figure}
\centering
\includegraphics[scale=0.15]{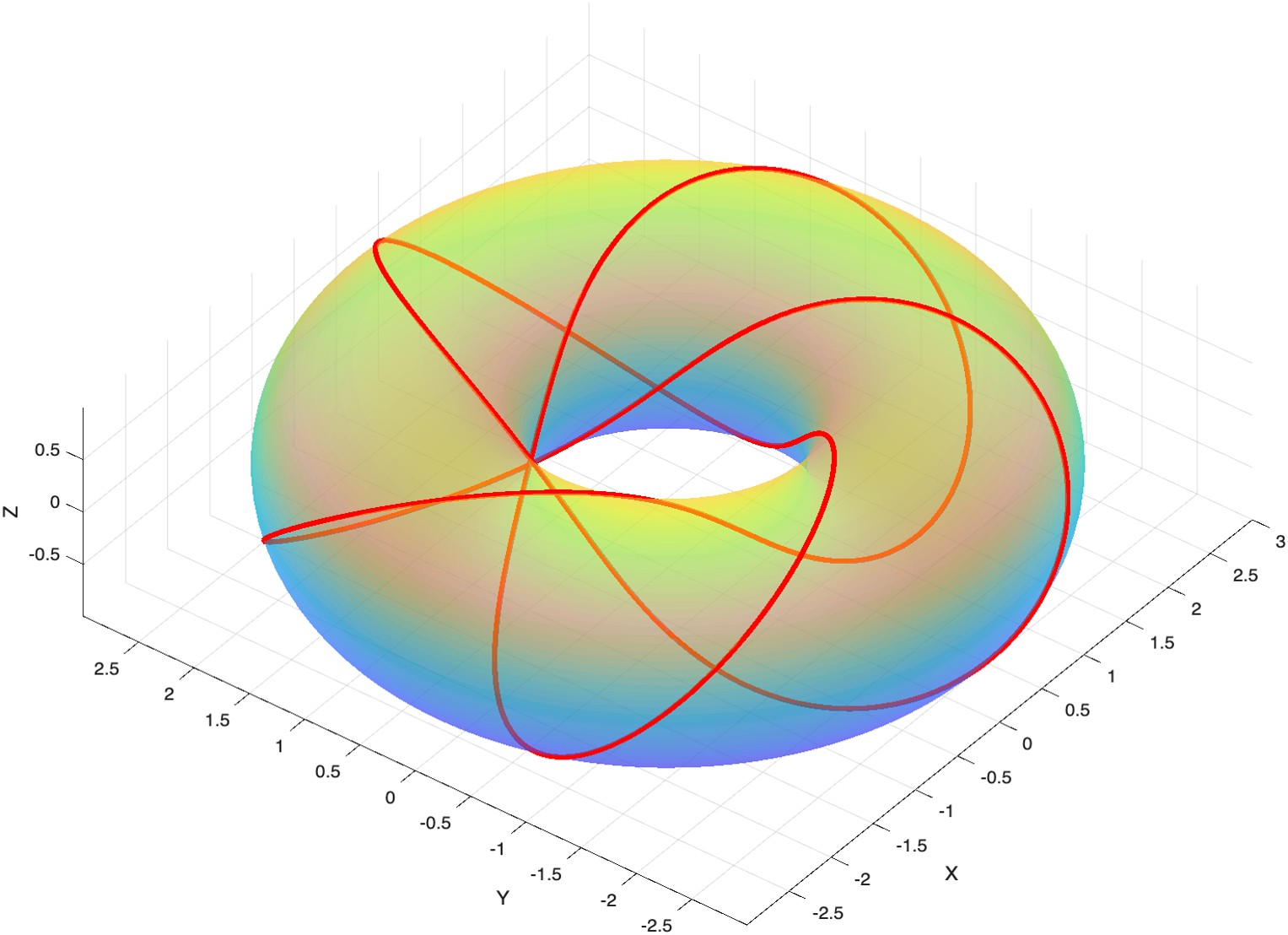}
\caption{\label{fig:toru} A closed loxodrome on a torus in Example~\ref{eg:torus} with $c=2$, corresponding to $\Omega = \frac{3}{5} \cdot 2\pi$ and azimuth $\theta = \arctan(\frac{3}{5}\sqrt{3})$.}
\end{figure}

The behavior of nonparallel maximal loxodromic curves is determined by the ratio $\Omega/(2\pi)$.
\begin{enumerate}
    \item If $\Omega$ is a rational multiple of $2\pi$, the loxodrome is a closed periodic curve. Moreover, when $\Omega = \frac{m}{n} 2\pi$ for some coprime integers $m, n$ with $n \not= 0$,  the loxodrome cycles through the profile curve $n$ times and completes $m$ full rotations around the rotation axis before returning to its starting point $p_0$. In Figure~\ref{fig:toru}, the loxodrome for $\Omega = \frac{3}{5} \cdot 2\pi$ cycles through the profile curve $5$ times and completes $3$ full rotations around the rotation axis.
    \item If $\Omega/2\pi$ is irrational, the loxodrome never closes. Instead, it loops endlessly around the torus, and its trajectory is dense in $S$.
\end{enumerate}

\begin{prop}\label{prop:closed_curve}
    Suppose $S$ is a torus-like surface. Every maximal parallel on $S$ is closed. A nonparallel maximal loxodromic curve $\lambda$ on $S$ is closed if and only if $\Omega/2\pi\in\mathbb Q$.
\end{prop}
\begin{proof}
    When the maximal loxodromic curve is a parallel, the result is clear. Let $\lambda(\tilde{v}) = \tilde\Phi(\tilde u(\tilde v),\tilde v)$ be a maximal loxodromic curve other than a parallel, where $\tilde u$ is given by (\ref{eq:u_of_v_tilde}). Suppose $\lambda$ is closed. Then $\lambda$ satisfies
    \begin{equation} \label{eq:pf_(prop:dense_curve)_1}
            \tilde{\Phi}(\tilde u(\tilde v),\tilde v) = \tilde{\Phi} (\tilde u(\tilde v+T),\tilde v+T),
    \end{equation}
    where $T$ is the period of $\lambda$. The third component of (\ref{eq:pf_(prop:dense_curve)_1}) yields $\tilde z(\tilde v)=\tilde z(\tilde v+T)$. Furthermore, squaring and adding the first two components yields $\tilde{r}(\tilde v)=\tilde{r}(\tilde v+T)$. Thus we have $\tilde{\alpha}(\tilde v) = \tilde{\alpha}(\tilde v+T)$. Since the closed profile curve is traversed exactly once on $[a,b]$, the periodic extension $\tilde{\alpha}$ has least period $L=b-a$. Hence $T=nL$ for some positive integer $n$. Applying (\ref{eq:quasi-periodic_of_u}), we have
    \begin{align*}
        \cos(\tilde u(\tilde v)) &= \cos(\tilde u(\tilde v+T)) = \cos(\tilde u(\tilde v)+n\Omega),\\
    \sin(\tilde u(\tilde v)) &= \sin(\tilde u(\tilde v+T)) = \sin(\tilde u(\tilde v)+n\Omega).
    \end{align*}
    By the periodicity of trigonometric functions, we have $ 2 \pi m = n\Omega$ for some integer $m$.

    Conversely, suppose $2 \pi m = n\Omega$. Reversing the preceding steps, let $T = n L$. We have $\tilde u(\tilde v+T)=\tilde u(\tilde v)+n\Omega = \tilde u(\tilde v)+2 \pi m$ by (\ref{eq:quasi-periodic_of_u}). Consequently, by the periodicity of the trigonometric functions and the profile curve, it follows that $\lambda(\tilde v) = \lambda(\tilde v+T)$.
\end{proof}
By Proposition~\ref{prop:closed_curve} and (\ref{eq:periodic_u}), a nonparallel maximal loxodromic curve is closed if its azimuth satisfies $\tan(\theta )= \frac{m}{n}\frac{2\pi}{\int_{a}^b \frac{1}{r(\tau)} \,d \tau}$.

Returning to the linearizing parametrization $\Psi$ described in \eqref{eq:Conformal_parametrization}, we can extend this coordinate system globally for torus-like surfaces. Since the function $r$ extends to a strictly positive periodic function $\tilde r$, the integral
\begin{equation}
    w(v) = \int_{v_0}^{v} \frac{1}{\tilde{r}(\tau)}\,d\tau
\end{equation}
is a strictly increasing diffeomorphism from $\mathbb{R}$ onto $\mathbb{R}$ for a fixed $v_0 \in I$. Let $v(w)$ be the inverse of $w$. We can define the extended linearizing parametrization $\tilde{\Psi}: \mathbb{R}^2 \to S$ by
\begin{equation}\label{eq:global_conformal_parametrization_on_torus}
    \tilde{\Psi}(u,w) = \tilde{\Phi}(u,v(w)).
\end{equation}
Since $\tilde{\Psi}$ is locally identical to $\Psi$, any loxodrome represented in the extended $(u,w)$-plane continues to satisfy the linear equations derived in Corollary~\ref{cor:uw}. Furthermore, the periodicity of $\tilde{\Phi}$ equips the $(u,w)$-plane with a lattice structure generated by the periods $2\pi$ in the $u$-direction and $\Omega_w = \int_a^b \frac{1}{r(\tau)} \, d\tau$ in the $w$-direction. In the $(u,w)$-plane, the points $(u,w)$ and $(u+2\pi m,w+n\Omega_w)$ represent the same point on $S$ for any integers $m,n$. Identifying these points gives a torus, and loxodromes on $S$ are represented by straight lines in this plane. By Lemma~\ref{lemma:maximal_solution_torus}, a nonparallel maximal loxodromic curve can be expressed as $\tilde{\Psi}(\tilde u(\tilde w),\tilde w)$, where
\begin{equation}\label{eq:uw_tilde}
    \tilde u(\tilde w) = \tan(\theta)(\tilde w - w_0) + \tilde u_0.
    \end{equation}
We now apply this global parametrization $\tilde{\Psi}$ to prove the following proposition.

\begin{prop}\label{prop:dense_loxo_on_torus}
    Suppose $S$ is a torus-like surface. Let $\lambda$ be a nonparallel maximal loxodromic curve on $S$. If $\Omega$ is not a rational multiple of $2\pi$, then $\lambda$ is not closed and its trajectory is dense in $S$.
\end{prop}
\begin{proof}
By the symmetry of the surface of revolution, we restrict our attention to $\theta \in (0,\pi/2)$. Let $F:\mathbb{R}^2\to \mathbb{T}^2$ be the map defined by $F(x,y)=(\exp(2\pi i x),\exp(2\pi i y))$. For a given $\alpha \not\in \mathbb{Q}$, define $G_\alpha:\mathbb{R}\to\mathbb{R}^2$ by $G_\alpha(t)=(t,\alpha t)$ and $H = F\circ G_\alpha$. According to Example 6.15 in \cite{boothby2003introduction}, the set
$$
D_\alpha=F^{-1}(H(\mathbb{R})) = \left\{(x,y):x=\frac{1}{\alpha}y+(n-\frac{1}{\alpha}m) \text{ for some } n, m \in \mathbb{Z}\right\}
$$
is dense in $\mathbb{R}^2$, which implies that $H(\mathbb{R})$ is dense in $\mathbb{T}^2$.

Without loss of generality, assume that the maximal loxodromic curve $\lambda$ passes through the point $\tilde{\Psi}(0,0)$, and that the longitude increment $\Omega$ is an irrational multiple of $2\pi$. We set $\alpha = 2\pi/\Omega$. Since $\Omega$ is an irrational multiple of $2\pi$, the value $\alpha$ is also irrational, which ensures that $D_\alpha$ is dense in $\mathbb{R}^2$. The set $\tilde{\Psi}^{-1}(\lambda(\mathbb{R}))$ can be described by
\begin{equation*}
    \tilde{\Psi}^{-1}(\lambda(\mathbb{R})) = \left\{(u,w) : u = \tan(\theta) w + (2\pi n - \Omega m) \text{ for some } n, m \in \mathbb{Z} \right\}.
\end{equation*}
Define the continuous mapping $f: \mathbb{R}^2 \to \mathbb{R}^2$ by
$$
f(x,y) = \left( 2\pi x , \Omega_w y \right).
$$
Then $f$ is a surjective continuous map on $\mathbb{R}^2$ with  $f(D_\alpha) = \tilde{\Psi}^{-1}(\lambda(\mathbb{R}))$. Since  $D_\alpha$ is dense in $\mathbb{R}^2$, it follows that the image $f(D_\alpha) = \tilde{\Psi}^{-1}(\lambda(\mathbb{R}))$ is dense in $\mathbb{R}^2$. Finally, since $\tilde{\Psi}$ is a continuous surjective mapping, we conclude that the trajectory $\lambda(\mathbb{R})$ of the maximal loxodromic curve $\lambda$ is dense in $S$.
\end{proof}

\begin{eg}\label{eg:torus}
    Consider a profile curve $\alpha(v) = (r(v),0,z(v))$, where $c>1$, $r(v)=\cos(v)+c$, and $z(v)=\sin(v)$, defined on the interval $I = [-\pi,\pi]$. By \eqref{eq:quasi-periodic_of_u}, the maximal loxodromic curve with azimuth $\theta \in (0,\pi/2)$ starting at $\Phi(0,0)$ is given by
\begin{equation}
    \tilde u(\tilde v) = \tan\theta \left( \mathcal{I}_c(\tilde v) + \frac{2\pi}{\sqrt{c^2-1}} \left\lfloor \frac{\tilde v+\pi}{2\pi} \right\rfloor \right),
\end{equation}
for each $\tilde v \in \mathbb R$, where $\mathcal{I}_c(\tilde v)$ is defined as
\begin{equation*}
    \mathcal{I}_c(\tilde v) =
    \begin{cases}
        \displaystyle \frac{2}{\sqrt{c^2-1}} \arctan \left( \sqrt{\frac{c-1}{c+1}} \tan\frac{\tilde v}{2} \right) & \text{if } \tilde v \neq (2k+1)\pi, \\
        \displaystyle -\frac{\pi}{\sqrt{c^2-1}} & \text{if } \tilde v = (2k+1)\pi,
    \end{cases}
\end{equation*}
    and the longitude increment is $\Omega = \frac{2\pi \tan(\theta)}{\sqrt{c^2-1}}$.
    If $\theta=\arctan\left(\frac{m}{n}\sqrt{c^2-1}\right)$ for some integers $m,n$ with $n\neq0$, then the loxodrome with azimuth $\theta$ is closed. Otherwise, its trajectory is dense in $S$. In Figure~\ref{fig:toru}, we set $c = 2$ and choose an azimuth $\theta = \arctan(\frac{3}{5}\sqrt{3})$. This loxodrome cycles through the profile curve $5$ times and completes $3$ full rotations around the rotation axis before returning to its starting point.
\end{eg}
\begin{eg}\label{eg:linearizing_parameterization_of_torus}
    Continuing from Example \ref{eg:torus}, we now construct the linearizing parametrization $\tilde{\Psi}$. Since $w(v) = \int_{0}^v \frac{1}{\tilde r(\tau)}\, d \tau$, where $\tilde r(v)=c+\cos v$ is the periodic extension of $r$ to $\mathbb R$, we have
    $$
    w(v) = \mathcal I_c(v) + \frac{2\pi}{\sqrt{c^2-1}} \left\lfloor\frac{v+\pi}{2\pi}\right\rfloor.
    $$
 For $v\neq(2k+1)\pi$, $k\in\mathbb Z$, we obtain $\tan\frac{\sqrt{c^2-1} w}{2} = \sqrt{\frac{c-1}{c+1}}\tan\frac v2$. Thus, we have $\tan\frac v2 = \sqrt{\frac{c+1}{c-1}} \tan\frac{\sqrt{c^2-1} w}{2}$. Using the half-angle identities,
    $
    \cos v=\frac{1-\tan^2(v/2)}{1+\tan^2(v/2)},
    $ and $
    \sin v=\frac{2\tan(v/2)}{1+\tan^2(v/2)},
    $
    we obtain
    $$
    \cos v(w)
    =
    \frac{c\cos(\sqrt{c^2-1} w)-1}
        {c-\cos(\sqrt{c^2-1} w)}, \quad \sin v(w)=\frac{\sqrt{c^2-1}\sin(\sqrt{c^2-1} w)}
     {c-\cos(\sqrt{c^2-1} w)}.
    $$
    Since $\tilde r(v)=c+\cos v$, we have
$$
\tilde r(v(w))
=
\frac{c^2-1}{c-\cos(\sqrt{c^2-1} w)}.
$$
Thus, the extended linearizing parametrization is

$$
\tilde{\Psi}(u,w ) = \left(
\frac{(c^2-1)\cos u}
{c-\cos\left(\sqrt{c^2-1}w\right)},
\frac{(c^2-1)\sin u}
{c-\cos\left(\sqrt{c^2-1}w\right)},
\frac{\sqrt{c^2-1}\sin\left(\sqrt{c^2-1}w\right)}
{c-\cos\left(\sqrt{c^2-1}w\right)}
\right).
$$
In addition, $\Omega_w = \int_a^b \frac{1}{r(\tau)} \, d\tau = \frac{2\pi}{\sqrt{c^2-1}}.$
\end{eg}

\subsection{Paraboloid-Like Surfaces}\label{subsec:paraboloidlike1}
A surface of revolution with only one pole and a non-self-intersecting profile curve is called a paraboloid-like surface. We now prove Lemma~\ref{lemma:modulo_azimuths} for this type of surface.
\begin{lemma}\label{lemma:unique_and_modulo_for_paraboloid}
    Suppose $S$ is a paraboloid-like surface with pole $p$. If $\lambda:J\to S$ is an arc-length-parametrized loxodromic curve passing through the pole $p$ with $\lambda(s_p)=p$ for some $s_p\in\operatorname{int}J$, then 
    \begin{enumerate}
        \item\label{item:unique}  The point $s_p$ is the unique point in $J$ whose image under $\lambda$ is $p$. 
        \item\label{item:modulo} The directed azimuths of all restrictions described in condition~\ref{item:loxodromic-local} of Definition~\ref{defi:loxodromic_curve} are congruent to $0$ modulo $\pi$.
    \end{enumerate}
\end{lemma}
\begin{proof}
 We first prove \eqref{item:unique}. Suppose that $\lambda(s_*)=p$ for some $s_*\neq s_p$.
Without loss of generality, suppose $s_*<s_p$.
We claim that there exists $s_1\in[s_*,s_p)$ such that
$\lambda(s_1)=p$ and $\lambda(s)\neq p$ for every $s\in(s_1,s_p)$.
Suppose otherwise. Starting with $t_0=s_*$, we can then choose
$t_{n+1}\in(t_n,s_p)$ with $\lambda(t_{n+1})=p$ for every $n\geq0$.
The increasing sequence $(t_n)$ is bounded above by $s_p$, so it
converges to some $c\in(s_*,s_p]\subseteq\operatorname{int}J$.
Since $\lambda$ is continuous at $c$, convergence $t_n\to c$ implies
$$
\lambda(c)=\lim_{n\to\infty}\lambda(t_n)=p.
$$
The definition of the derivative therefore gives
$$
\lambda'(c)
=\lim_{n\to\infty}\frac{\lambda(t_n)-\lambda(c)}{t_n-c}
=0,
$$
contradicting the regularity of $\lambda$ at $c$.
Thus such an $s_1$ exists, and
$\lambda|_{(s_1,s_p)}$ is a loxodrome. Let $\theta$ be its directed azimuth, and write $\lambda(s)=\Phi(u(s),v(s))$ for $s\in(s_1,s_p)$. By \eqref{eq:dvds}, we have $v'(s)=\cos\theta$, and hence $v(s)=\cos\theta\,s+C$ for some constant $C$. Let $v_p\in I$ satisfy $p=\Phi(0,v_p)$.
Since $\lambda$ is continuous at $s_1$ and $s_p$, we obtain
$$
\lim_{s\to s_1^+}\lambda(s)
=
\lim_{s\to s_p^-}\lambda(s)=p.
$$
Using $\lambda(s)=\Phi(u(s),v(s))$ and
$p=\alpha(v_p)=(0,0,z(v_p))$, we obtain
$$
\|\lambda(s)-p\|^2
=
r(v(s))^2+\bigl(z(v(s))-z(v_p)\bigr)^2
=
\|\alpha(v(s))-\alpha(v_p)\|^2.
$$
Thus $\lambda(s)\to p$ implies
$\alpha(v(s))\to\alpha(v_p)$ at both endpoints.
Since $S$ is regular at $p$ and the profile has no
self-intersections, profile points approaching $\alpha(v_p)$
have parameters approaching $v_p$. Hence
$$
\lim_{s\to s_1^+}v(s)
=
\lim_{s\to s_p^-}v(s)=v_p.
$$
Substituting $v(s)=\cos\theta\,s+C$, we obtain
$$
\cos\theta\,s_1+C
=
\cos\theta\,s_p+C
=
v_p.
$$
Since $s_1<s_p$, it follows that $\cos\theta=0$ and
$v(s)\equiv v_p$. Hence
$\lambda(s)=\Phi(u(s),v_p)=p$ for every $s\in(s_1,s_p)$,
contradicting the choice of $s_1$. Therefore, $s_p$ is the unique point in $J$ mapped to $p$ by $\lambda$.

We next prove \eqref{item:modulo}. First consider $J_0=(s_0,s_p)$, where $s_0\in\operatorname{int}J$ and $s_0<s_p$. By \eqref{item:unique}, $\lambda|_{J_0}$ is a loxodrome. Let $\theta$ be its directed azimuth. Suppose, for a contradiction, that
$\theta\not\equiv0 \bmod\pi$. We can write $\lambda|_{J_0}=\Phi(u(s),v(s))$ for $s\in J_0$, where $u$ and $v$ satisfy \eqref{eq:duds} and \eqref{eq:dvds}. Suppose the pole corresponds to the endpoint $a$ of the interval of the profile curve. Since $\lambda(s)\to p$ as $s\to s_p^-$, we have $v(s)\to a$. We must have $\cos\theta\neq0$, since otherwise $v$ would be constant and $\lambda(s)$ could not approach $p$. Since $\theta$ is not congruent to $0$ modulo $\pi$, we also have $\sin\theta\neq0$. Hence $\frac{du}{dv}=\frac{\tan\theta}{r(v)}$. Since $a$ is the left endpoint of $I$, the conditions $r(v)\geq0$, $r(a)=0$, and $|r'(a)|=1$ imply that $r'(a)=1$. We have
$\lim_{v\to a^+}r(v)/(v-a)=1$.
Therefore, there exists $\delta>0$ such that
$0<r(v)<2(v-a)$ for every $v\in(a,a+\delta)$. It follows that
$\int_v^{a+\delta}1/r(\tau)\,d\tau>
\int_v^{a+\delta}1/[2(\tau-a)]\,d\tau$,
which tends to infinity as $v\to a^+$. Since $\tan\theta\neq0$, we conclude that $u(v)$ diverges to $\pm\infty$ as $v\to a^+$. Thus $u(s)\to+\infty$ or $u(s)\to-\infty$
as $s\to s_p^-$. Moreover,
\begin{align*}
\lambda'(s)=
\begin{pmatrix}
-\sin\theta\sin u(s)+\cos\theta\,r'(v(s))\cos u(s)\\
\sin\theta\cos u(s)+\cos\theta\,r'(v(s))\sin u(s)\\
\cos\theta\,z'(v(s))
\end{pmatrix}.
\end{align*}
Suppose first that $u(s)\to+\infty$ as $s\to s_p^-$. Since $u(s)$ is continuous, strictly monotonic, and unbounded, there exists an integer $N$ such that, for every integer $n\geq N$, there exist $s_n,t_n\in J_0$ satisfying $u(s_n)=2n\pi$ and $u(t_n)=2n\pi+\pi$, and both sequences $(s_n)$ and $(t_n)$ converge to $s_p$. The case in which $u(s)\to-\infty$ is similar. Since $r'(v(s))\to1$ and $z'(v(s))\to0$, $$
\lim_{n\to\infty}\lambda'(s_n)
=(\cos\theta,\sin\theta,0),
\quad
\lim_{n\to\infty}\lambda'(t_n)
=(-\cos\theta,-\sin\theta,0).
$$
These limits are distinct because the two vectors are nonzero and are negatives of one another. Hence $\lambda'(s)$ does not converge as $s\to s_p^-$. This contradicts the smoothness of $\lambda$ at $s_p$. Therefore, $\theta$ is congruent to $0$ modulo $\pi$. The case where the pole corresponds to $b$ is analogous. Reversing the parameter gives the same conclusion for intervals of the form $(s_p,s_0)$.

Now let $\lambda|_{J_0}$ be any other restriction, and choose $s_0\in J_0$. Without loss of generality, suppose $s_0<s_p$. By \eqref{item:unique}, we have $\lambda(s)\neq p$ for every $s\in(s_0,s_p)$. Hence $\lambda|_{(s_0,s_p)}$ is a loxodrome whose directed azimuth is congruent to $0$ modulo $\pi$. Since $(s_0,s_p)\cap J_0\neq\emptyset$, the restriction $\lambda|_{J_0}$ has the same directed azimuth. Therefore, its directed azimuth is also congruent to $0$ modulo $\pi$.
\end{proof}

\begin{proof}[Proof of Lemma~\ref{lemma:modulo_azimuths} for paraboloid-like surfaces]
Let $\lambda:J\to S$ be a loxodromic curve. If $\lambda$ does not pass through $p$ in the interior of its domain, the proof is the same as that for cylinder-like surfaces.

If $\lambda$ passes through $p$, we parametrize $\lambda$ by arc length. By Lemma~\ref{lemma:unique_and_modulo_for_paraboloid}, the directed azimuths of all restrictions described in condition~\ref{item:loxodromic-local} are congruent to $0$ modulo $\pi$. Since a regular reparametrization does not change the directed azimuth modulo $\pi$, the same conclusion holds for the original curve.
\end{proof}
If the azimuth of a maximal loxodromic curve is $\pi/2$ modulo $\pi$, then \eqref{eq:dvds} gives $dv/ds=0$. Hence $v$ is constant, and the curve is contained in a parallel. Since the solution of \eqref{eq:duds} is defined for every $s\in\mathbb R$, maximality implies that the curve is the entire parallel. Therefore, it is closed.

Let $\lambda:J\to S$ be a maximal loxodromic curve whose azimuth is neither $0$ nor $\pi/2$ modulo $\pi$. By Lemma~\ref{lemma:unique_and_modulo_for_paraboloid}, $\lambda$ does not pass through the pole. On $\operatorname{int}J$, write $\lambda(s)=\Phi(u(s),v(s))$, and let $\theta$ be the directed azimuth of this restriction. Since $\cos\theta\neq0$, the function $v(s)$ is strictly monotonic. We may therefore use $v$ instead of $s$ and write $\lambda(v)=\Phi(u(v),v)$. By Corollary~\ref{cor:uv},
$u(v)=\tan\theta\int_{v_0}^{v}1/r(\tau)\,d\tau+u_0$.
By maximality, $v$ ranges over $\operatorname{int}I$.

Suppose that the pole corresponds to the endpoint $a$ of $I$. As $v\to a$, the curve approaches the pole $p$, and the proof of Lemma~\ref{lemma:unique_and_modulo_for_paraboloid} gives $|u(v)|\to\infty$. Hence $\lambda$ winds around the axis of rotation infinitely many times. If $\lambda$ did not reach $p$, adjoining $p$ as an endpoint would give a loxodromic curve on a larger interval, contradicting maximality. Therefore, $\lambda$ has $p$ as an endpoint. Finally, $\lambda$ is non-self-intersecting. Indeed, suppose that $\lambda(s_1)=\lambda(s_2)$ for some $s_1,s_2\in\operatorname{int}J$. Since the profile curve is non-self-intersecting, we have $v(s_1)=v(s_2)$. Since $v$ is strictly monotonic, it follows that $s_1=s_2$. Moreover, no point of $\operatorname{int}J$ is mapped to $p$, because $\lambda$ does not pass through the pole. Hence $\lambda$ is non-self-intersecting on $J$.

Let $\lambda:J\to S$ be a maximal loxodromic curve whose azimuth is $0$ modulo $\pi$. Away from the pole, \eqref{eq:duds} gives $du/ds=0$, so $\lambda$ follows a fixed-longitude meridian. Maximality forces this meridian to continue until it reaches the pole. Since a meridian continues smoothly through the pole along the opposite longitude, maximality also forces $\lambda$ to contain this continuation. Hence the pole is the image of an interior point of $J$, and therefore $\lambda$ passes through it.

Without loss of generality, suppose that the pole corresponds to $v=b$. Since the azimuth of $\lambda$ is $0$ modulo $\pi$, \eqref{eq:dvds} gives $dv/ds=\pm1$. Hence $v(s)=\pm s+C$ before $\lambda$ reaches the pole. Changing the sign of $s$ and adding a constant if necessary, we may write $v(s)=s$. Therefore,
\begin{equation}\label{eq:meridian_for_1pole}
\lambda(s)=
\begin{cases}
\Phi(u_0,s), & s\in I, \\
\Phi(u_0+\pi,2b-s), & 2b-s\in I\setminus\{b\}.
\end{cases}
\end{equation}
Here, a maximal loxodromic curve with azimuth $0$ modulo $\pi$ is called a \textbf{maximal meridian}.

\begin{figure}
\centering
\includegraphics[scale=0.4]{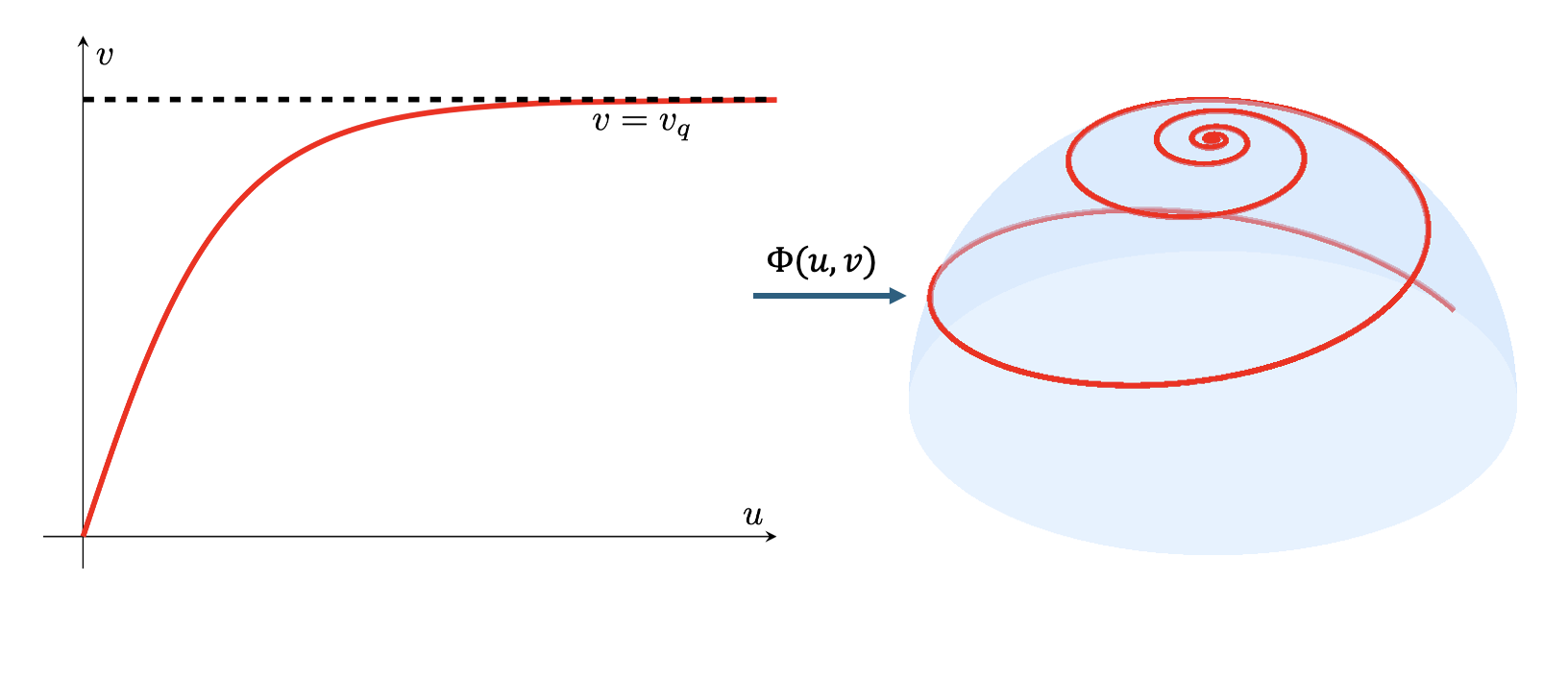}
\caption{\label{fig:hemisphere}
A nonmeridional maximal loxodrome approaching the pole
of a hemisphere. Adjoining the pole as an endpoint gives
a maximal loxodromic curve. The left figure shows the
loxodrome in the $(u,v)$-plane, and the right figure shows
its image under $\Phi$.}
\end{figure}

\begin{cor}\label{cor:topo_of_loxodrome_on_PS}
    Let $S$ be a paraboloid-like surface with a pole $q$. There are only three kinds of maximal loxodromic curves:
    \begin{enumerate}
    \item closed maximal parallels,
    \item maximal meridians passing smoothly through $q$, and
    \item non-self-intersecting maximal loxodromic curves that reach $q$ but do not pass through it.
\end{enumerate}
\end{cor}

\begin{eg}\label{eg:hemisphere}
    Given a profile curve $\alpha(v) = (\cos v, 0, \sin v)$ on the interval $I=[0,\pi/2]$, a parametrization of the closed upper hemisphere $S$ is given by
    $$
    \Phi(u,v) = (\cos (v) \cos (u), \cos (v) \sin (u), \sin (v)).
    $$
    By Corollary~\ref{cor:uv}, the maximal loxodromic curve that passes through $p_0 = \Phi(u_0,v_0)$ and has azimuth neither $0$ nor $\pi/2$ modulo $\pi$ is given by $\lambda(v) = \Phi(u(v),v)$, where
    $$
    u(v) = \tan(\theta) \big( \ln |\sec (v) + \tan (v)| - \ln |\sec (v_0) + \tan (v_0)| \big) + u_0
    $$
    for $0\leq v<\pi/2$, with the north pole included
    as an endpoint of the maximal loxodromic curve. The arc length $L$ of the segment of $\lambda$ from $v_0$ to $v$ on the unit hemisphere is given by $L = |\sec \theta| \, |v - v_0|$. As $v\to\frac{\pi}{2}^{-}$, we have $|u(v)|\to\infty$. Consequently, the maximal loxodromic curve winds around the \(z\)-axis infinitely many times while converging to the north pole. Its total length from \(v=v_0\) to the pole is $|\sec\theta| \left(\frac{\pi}{2}-v_0\right)<\infty$. Hence, the north pole is a finite-length endpoint of the maximal loxodromic curve. By Lemma~\ref{lemma:unique_and_modulo_for_paraboloid}, the curve cannot pass smoothly through the north pole, as illustrated in Figure~\ref{fig:hemisphere}.
\end{eg}
\subsection{Sphere-Like Surfaces}\label{subsec:spherelike1}
A sphere-like surface is a surface whose profile curve intersects the axis of rotation at both endpoints. According to equation (\ref{eq:G-B_Thm}), for a sphere-like surface, the regularity at the poles and the assumption of positive orientation imply $r'(a)=1$ and  $r'(b)=-1$. Thus, $\int_S K\,dA=4\pi$. Since $S$ is a compact, connected, orientable surface without boundary, the Gauss--Bonnet theorem
gives $\chi(S)=2$. By the classification of compact connected
orientable surfaces, $S$ is homeomorphic to a sphere
$\mathbb{S}^2$. Consequently, its fundamental group is trivial.

The argument used in the proof of Lemma~\ref{lemma:unique_and_modulo_for_paraboloid} applies near each pole. Thus, a maximal loxodromic curve can pass through a pole only if its azimuth is $0$ modulo $\pi$, and the proof of Lemma~\ref{lemma:modulo_azimuths} for sphere-like surfaces is the same as that for paraboloid-like surfaces.

Now let $\lambda:J\to S$ be a maximal loxodromic curve whose azimuth is neither $0$ nor $\pi/2$ modulo $\pi$. On $\operatorname{int}J$, let $\theta$ be its directed azimuth. Since $\cos\theta\neq0$, \eqref{eq:dvds} shows that $v$ is strictly monotonic. Maximality forces $v$ to range over $\operatorname{int}I$; otherwise, the solution of \eqref{eq:duds} and \eqref{eq:dvds} would give a loxodromic curve on a larger interval. Hence $\lambda$ approaches one pole as $v\to a$ and the other as $v\to b$. Both poles belong to $\lambda$. Otherwise, adjoining either pole as an endpoint would give a loxodromic curve on a larger interval, contradicting maximality. Thus, the two poles are the endpoints of $\lambda$. Moreover, if $\lambda(s_1)=\lambda(s_2)$ for $s_1,s_2\in\operatorname{int}J$, the non-self-intersection of the profile curve gives $v(s_1)=v(s_2)$. Since $v$ is strictly monotonic, we have $s_1=s_2$. Therefore, $\lambda$ is non-self-intersecting.

\begin{figure}
\centering
\includegraphics[scale=0.15]{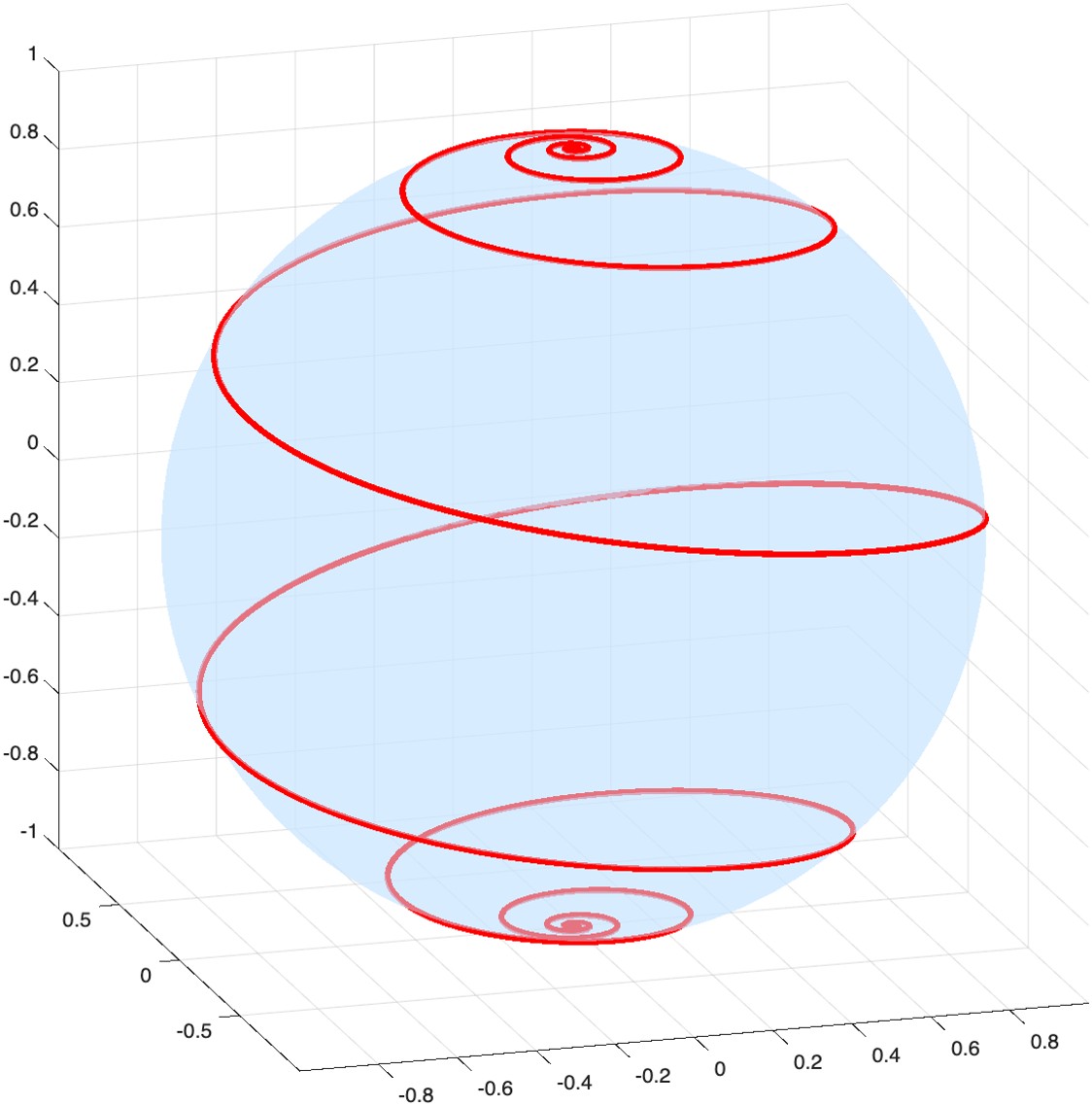}
\caption{\label{fig:sphere} A maximal loxodromic curve on a sphere.}
\end{figure}

Let $L = b-a$ be the arc length of the profile curve between the two poles. By applying the smooth continuation described in the preceding subsection at both poles, we can extend each fixed-longitude meridian to a maximal meridian. After passing through a pole, the continuation follows the opposite longitude. After passing through both poles, it returns to its original longitude and starting point. Thus, every maximal meridian is a smooth closed loxodromic curve of period $2L$ and azimuth $0$. We can write a maximal meridian with an initial longitude $u_0$ as
\begin{equation}\label{eq:meridian_for_2pole}
    \lambda(v) = \begin{cases}
        \Phi(u_0, v-2nL), & \text{if } v-2nL \in I, \text{ for some } n \in \mathbb{Z}, \\
        \Phi(u_0 + \pi, 2b - v+2nL), & \text{if } 2b - v+2nL \in I \setminus \{b\}, \text{ for some } n \in \mathbb{Z}.
    \end{cases}
\end{equation}
\begin{cor}\label{cor:topo_of_loxodrome_on_SS}
    Let $S$ be a sphere-like surface. There are only three kinds of maximal loxodromic curves:
    \begin{enumerate}
        \item closed maximal parallels,
        \item maximal meridians passing smoothly through the two poles, and
        \item non-self-intersecting maximal loxodromic curves having the two poles as their endpoints.
    \end{enumerate}
\end{cor}
\begin{eg}\label{eg:sphere}
    Consider the sphere generated by the profile curve $\alpha(v)=(\cos v,0,\sin v)$ on $I=[-\pi/2,\pi/2]$. A maximal loxodromic curve that is neither a meridian nor a parallel is described away from the poles by the same formula as in Example~\ref{eg:hemisphere}. As $v\to\pm\pi/2$, we have $|u(v)|\to\infty$, and the maximal loxodromic curve has the north and south poles as its endpoints, as illustrated in Figure~\ref{fig:sphere}.
\end{eg}
\section{The Direct and Inverse Problems}\label{sec:direct_inverse}
For each admissible choice of coordinate representatives considered below, the corresponding formula determines a loxodromic segment joining the prescribed points. An \textbf{(oriented) loxodromic segment} from an initial point $p_0$ to a terminal point $p_1$ is an oriented restriction $\lambda: [s_0,s_1] \to S$ of a maximal loxodromic curve, where $s_0<s_1$, $\lambda(s_0)=p_0$ and $\lambda(s_1)=p_1$. The azimuth of the oriented loxodromic segment is its initial directed azimuth, defined as the directed azimuth of $\lambda$ on $(s_0, s_0+\varepsilon)$, where $\varepsilon>0$ is sufficiently small that no point of this interval is mapped to a pole. The \textbf{loxodromic distance} from $p_0$ to $p_1$ is the minimum
of the lengths of all loxodromic segments from $p_0$ to $p_1$. We use these loxodromic segments to address the direct and inverse problems described in the Introduction.  

\subsection{The Direct Problem}\label{directproblem1}
Consider a loxodromic segment $\lambda$ of directed azimuth $\theta$ and length $d$, starting at the point $p_0=\Phi(u_0,v_0)$. We assume that the initial point $p_0$ of a loxodromic segment is not a pole. The terminal point $p_1$, however, may be a pole. Consequently, the directed azimuth of the segment is well defined as its directed azimuth at $p_0$. We seek to determine its terminal point $p_1$ along with the range of valid lengths $d$, depending on the type of the surface $S$.

\subsubsection*{(a)}
Suppose $S$ is a cylinder-like surface. By Lemma~\ref{lemma:maximal_solution_represents_loxodromic_curve}, we can write $\lambda(s)=\Phi(u(s),v(s))$, where $u$ and $v$ satisfy \eqref{eq:duds} and \eqref{eq:dvds}. Write $p_1=\Phi(u_1,v_1)$. Since $s$ is the arc-length parameter, we have $d=s_1-s_0$. Integrating \eqref{eq:dvds} gives $v_1=d\cos\theta+v_0$. If $\theta = \pi/2$ or $3\pi/2$, we have
\begin{align}\label{eq:second_point_for_parl}
    v_1 = v_0 ,\quad u_1 = u_0 + \frac{d\sin\theta}{r(v_0)}.
\end{align}
Otherwise, the coordinates of the terminal point are given by
\begin{align}
    u_1 &= \int_{v_0}^{d\cos\theta + v_0} \frac{\tan\theta}{r(\tau)} \,d\tau + u_0, \label{eq:second_point_for_u} \\
    v_1 &= d\cos\theta + v_0. \label{eq:second_point_for_v}
\end{align}
In this case, the length $d$ is bounded by the endpoints of the interval $I$, which is the domain of the profile curve. Specifically, when the curve travels in the direction of increasing $v$ ($\cos\theta > 0$), we must have $d \leq (b-v_0)/\cos\theta$ (with strict inequality if $S$ does not contain the upper boundary). Similarly, when the curve travels in the direction of decreasing $v$ ($\cos\theta < 0$), we have $d \leq (a-v_0)/\cos\theta$ (with strict inequality if $S$ does not contain the lower boundary). For parallels, the length $d$ is unbounded.

\subsubsection*{(b)}
For a torus-like surface, the terminal point is determined by the same formulas as in the cylinder-like case, with $r$ replaced by the periodic function $\tilde r$. 
Write $p_1=\tilde\Phi(u_1,v_1)$.
If the segment is contained in a parallel, then
$\theta=\pi/2$ or $3\pi/2$, and the coordinates
of the terminal point are
\begin{align}\label{eq:second_point_for_parl_torus}
    v_1 = v_0 ,\quad u_1 = u_0 + \frac{d\sin\theta}{\tilde r(v_0)}.
\end{align}
Otherwise, the coordinates of the terminal point are
\begin{align}
    u_1 &= \int_{v_0}^{d\cos\theta + v_0} \frac{\tan\theta}{\tilde r(\tau)} \,d\tau + u_0, \label{eq:second_point_for_u_torus} \\
    v_1 &= d\cos\theta + v_0. \label{eq:second_point_for_v_torus}
\end{align}
The value $v_0+d\cos\theta$ is first computed on $\mathbb R$, and its value modulo $L$ determines the profile coordinate of the terminal point. Since every maximal loxodromic curve on $S$ is defined on $\mathbb R$, the length $d$ has no upper bound for any directed azimuth $\theta$.

\subsubsection*{(c)}
Suppose that $S$ is a paraboloid-like surface with a pole $q=\Phi(u_q,b)$. If $\theta\neq0$ and the terminal point $p_1$ is not the pole, the formulas for $p_1$ and the bounds on $d$ are the same as in the cylinder-like case. In the special case where $\theta\neq0$, $\cos\theta>0$, and $d = \frac{b-v_0}{\cos\theta},$ the loxodromic segment terminates at the pole, and we define $p_1=q$. In this case, $v_1=b$, but $u_1$ is undefined since the corresponding integral diverges.

For $\theta=0$, the loxodromic segment is contained in a maximal meridian that passes through the pole. By \eqref{eq:meridian_for_1pole}, the terminal point is given by
\begin{align}\label{eq:second_point_for_v_para}
p_1=
\begin{cases}
\Phi(u_0, v_0+d), & \text{if } v_0+d \in I, \\
\Phi(u_0+\pi, 2b-(v_0+d)), & \text{if } 2b-(v_0+d) \in I\setminus\{b\}.
\end{cases}
\end{align}
If both endpoints of $I$ are included in the domain of the profile curve, the upper bound for $d$ is $2b-a-v_0$. If the endpoint corresponding to $v=a$ is excluded, then $d < 2b-a-v_0$. The case where $v_q=a$ is analogous.
\subsubsection*{(d)}
Suppose that $S$ is a sphere-like surface with poles $q_a=\Phi(u,a)$ and $q_b=\Phi(u,b)$. If a nonmeridional loxodromic segment does not terminate at a pole, its endpoint is determined by the same formulas and admissible length restrictions as in the cylinder-like case. If it terminates at a pole, its endpoint must be defined separately. Specifically, if $\cos\theta>0$ and
$
d=\frac{b-v_0}{\cos\theta},
$
then $p_1=q_b$. If $\cos\theta<0$ and
$
d=\frac{a-v_0}{\cos\theta},
$
then $p_1=q_a$. In either case, the longitude $u_1$ is undefined at the pole. For loxodromic segments contained in maximal meridians, \eqref{eq:meridian_for_2pole} gives the following terminal point when $\theta = 0$:
\begin{align}\label{eq:second_point_for_v_sphere_t}
    p_1 = \begin{cases}
        \Phi(u_0, (v_0+d)-2nL), & \text{if } (v_0+d)-2nL \in I \text{ for some } n \in \mathbb{Z}, \\
        \Phi(u_0 + \pi, 2b - (v_0+d)+2nL), & \text{if } 2b - (v_0+d)+2nL \in I \setminus \{b\} \text{ for some } n \in \mathbb{Z},
    \end{cases}
\end{align}
and for $\theta = \pi$ we have
\begin{align}\label{eq:second_point_for_v_sphere_b}
    p_1 = \begin{cases}
        \Phi(u_0, (v_0-d)-2nL), & \text{if } (v_0-d)-2nL \in I \text{ for some } n \in \mathbb{Z}, \\
        \Phi(u_0 + \pi, 2b - (v_0-d)+2nL), & \text{if } 2b - (v_0-d)+2nL \in I \setminus \{b\} \text{ for some } n \in \mathbb{Z}.
    \end{cases}
\end{align}
For both $\theta=0$ and $\theta=\pi$, the length $d$ is unbounded. Table~\ref{tab:direct_inverse_summary_refs} summarizes the formulas for the direct problem, with the conditions on the coordinates and the bounds on $d$ given above.

\subsection{The Inverse Problem}\label{inverseproblem1}
Given two distinct points $p_0, p_1$ on a surface of revolution, the inverse problem is to determine the loxodromic segments joining them, their possible directed azimuths, and the minimum length among those segments. For a loxodromic segment joining two nonpole points $p_0,p_1$ without passing through a pole, choose continuous coordinate functions $(u(s),v(s))$ along the segment, with endpoint values $(u_0,v_0)$ and $(u_1,v_1)$, respectively. By \eqref{eq:duds} and \eqref{eq:dvds}, we have
$$
u_1-u_0=\sin\theta\int_{s_0}^{s_1}
\frac{1}{r(v(s))}\,ds,
\quad
v_1-v_0=(s_1-s_0)\cos\theta.
$$
Since $r(v(s))>0$ and $s_1>s_0$, it follows that
$$
\operatorname{sgn}(\sin\theta)
=\operatorname{sgn}(u_1-u_0),
\quad
\operatorname{sgn}(\cos\theta)
=\operatorname{sgn}(v_1-v_0),
$$
where $\operatorname{sgn}$ is the sign function.

\subsubsection*{(a)}
Let $S$ be a cylinder-like surface. If $v_0=v_1$, any loxodromic segment joining $p_0$ to $p_1$ is contained in a maximal parallel, with directed azimuth $\theta=\pi/2$ or $3\pi/2$. By \eqref{eq:duds} and \eqref{eq:dvds}, the segment admits the arc-length parametrization $\lambda(s)=\Phi(u(s),v_0)$ for $s\in[s_0,s_1]$, where $$u(s)=u_0\pm\frac{s-s_0}{r(v_0)}.$$ Since $\lambda(s_1)=p_1$, its length is
$$
    d=s_1-s_0=r(v_0)|u_1-u_0|,
$$
where $u_0$ and $u_1$ are longitude representatives
chosen consistently with the segment.

Since longitude is defined modulo $2\pi$, different
windings around the axis give the possible lengths
\begin{equation}\label{eq:distance_of_loxodrome1}
d_n=r(v_0)|u_1-u_0+2\pi n|
\end{equation}
for each $n\in\mathbb Z$, where $u_0$ and $u_1$ are now fixed longitude
representatives of $p_0$ and $p_1$, respectively,
and only values with $d_n>0$ are considered. Consequently, the loxodromic distance is
$$
d=\min_{n\in\mathbb Z}
r(v_0)|u_1-u_0+2\pi n|.
$$

Assume $v_0 \neq v_1$. According to Corollary~\ref{cor:uv}, the parametrization of the loxodromic segment from $p_0$ to $p_1$ is given by $\Phi(u(v),v)$, where $u(v)$ satisfies \eqref{eq:dudv}. Thus,
\begin{equation}\label{eq:inverse_theta}
    u_1 - u_0 = \tan\theta \int_{v_0}^{v_1} \frac{1}{r(\tau)} \,d\tau.
\end{equation}
Rearranging equation (\ref{eq:inverse_theta}), we obtain
\begin{equation*}
    \tan\theta = \frac{u_1 - u_0}{\int_{v_0}^{v_1} \frac{1}{r(\tau)} \,d\tau}.
\end{equation*}
Also, we can parametrize the loxodromic segment $\lambda$ by $\Phi(u(s),v(s))$, where $u(s)$ and $v(s)$ satisfy \eqref{eq:duds} and \eqref{eq:dvds}. In particular, \eqref{eq:dvds} gives $v_1 - v_0 = d\cos\theta$. The arc length of $\lambda$ is given by
\begin{equation*}
d = |\sec\theta| |v_1 - v_0| = |v_1-v_0|\sqrt{1+(u_1 - u_0)^2\left({\int_{v_0}^{v_1} \frac{1}{r(\tau)} \,d\tau}\right)^{-2}}.
\end{equation*}
Similarly, the arc lengths depend on the choice of the coordinate $u$, and, for each integer $n\in\mathbb{Z}$, the corresponding azimuth satisfies 
\begin{equation}\label{eq:tangent_of_the_azimuth}
\tan\theta = \frac{u_1-u_0+2\pi n}
{\int_{v_0}^{v_1}
\frac{1}{ r(\tau)}\,d\tau}.
\end{equation}
For fixed $u_0$ and $u_1$, all possible arc lengths are
\begin{equation}\label{eq:distance_of_loxodrome2}
d_n = |v_1-v_0|\sqrt{1+(u_1 - u_0+2\pi n)^2\left({\int_{v_0}^{v_1} \frac{1}{r(\tau)} \,d\tau}\right)^{-2}}
\end{equation}
for each integer $n$. Therefore, the loxodromic distance is obtained by choosing $n$ to minimize $|u_1-u_0+2\pi n|$. Such an integer exists because every real number has an integer nearest to it.

Every loxodromic segment is either parallel or nonparallel, so the preceding cases exhaust all possibilities.

\subsubsection*{(b)}
Let $S$ be a torus-like surface. Since the profile curve is periodic with period $L$, the coordinate differences may be represented as $v_1-v_0+mL$ and $u_1-u_0+2\pi n$, respectively, where $m$ and $n$ are integers. Consequently, the possible azimuths and arc lengths of loxodromic segments from $p_0$ to $p_1$ depend on both coordinate differences. For fixed $u_0$ and $u_1$, the set of all possible lengths of nonparallel loxodromic segments is
\begin{equation}\label{eq:distance_of_loxodrome3}
d_{n,m} = |v_1 - v_0 + mL| \sqrt{1 + (u_1 - u_0 + 2\pi n)^2 \left( \int_{v_0}^{v_1+mL} \frac{1}{\tilde{r}(\tau)} \, d\tau \right)^{-2}}
\end{equation}
for all integers $m$ and $n$ satisfying $v_1-v_0+mL\neq0$. The minimum over all such pairs $(m,n)$ is attained
since the lengths tend to infinity as $|m|+|n|\to\infty$. For each such pair $(m,n)$, the corresponding azimuth satisfies
\begin{equation}\label{eq:tangent_of_the_azimuth_torus}
    \tan\theta=\frac{u_1-u_0+2\pi n}{\int_{v_0}^{v_1+mL}\frac{1}{\tilde r(\tau)}\,d\tau}.
\end{equation}
The directed azimuth $\theta\in[0,2\pi)$ is selected by
$\operatorname{sgn}(\cos\theta)
=\operatorname{sgn}(v_1-v_0+mL)$.

If $v_0\equiv v_1\pmod L$, there are also parallel segments with azimuth $\pi/2$ or $3\pi/2$. The loxodromic distance is the minimum among the lengths given by \eqref{eq:distance_of_loxodrome1} and \eqref{eq:distance_of_loxodrome3}. If $v_0\not\equiv v_1\pmod L$, no parallel segment joins the two points, and the loxodromic distance is the minimum among the lengths in \eqref{eq:distance_of_loxodrome3}.

Every loxodromic segment is either parallel or nonparallel, and the integers $m$ and $n$ describe all possible coordinate differences. Hence the preceding cases exhaust all possibilities.

\subsubsection*{(c)}
Suppose that $S$ is a paraboloid-like surface with pole $q$
corresponding to $v=v_q$. Suppose first that $p_0$ and $p_1$ are distinct nonpole points. Parallel segments and other nonmeridional loxodromic segments
are determined as in the cylinder-like case. Let $d_{N}$ denote the length of the shortest nonmeridional loxodromic segment joining $p_0$ and $p_1$. Thus, $d_{N}$ is determined by
\eqref{eq:distance_of_loxodrome1} when $v_0=v_1$, and by
\eqref{eq:distance_of_loxodrome2} when $v_0\neq v_1$.
In the latter case, the minimum is taken only over
integers $n$ satisfying $u_1-u_0+2\pi n\neq0$.

Meridian segments from $p_0$ to $p_1$ occur when
$u_1-u_0\equiv0\pmod\pi$. We distinguish two cases.

First, suppose that $u_1-u_0\equiv0\pmod{2\pi}$. Then $p_0$ and $p_1$ lie on the same fixed-longitude meridian branch. The meridian segment joining them has length
\begin{equation}\label{eq:distance_of_same_branch_meridian}
d_0=|v_1-v_0|.
\end{equation}
Its initial directed azimuth is $0$ when $v_1>v_0$, and $\pi$ when
$v_1<v_0$. In this case, the loxodromic distance is
$\min\{d_0,d_N\}$.

Next, suppose that $u_1-u_0\equiv\pi\pmod{2\pi}$. Then the meridian
segment joining $p_0$ and $p_1$ passes through the pole $q$. Its length
is
\begin{equation}\label{eq:distance_of_meridian}
d_q=|v_q-v_0|+|v_q-v_1|.
\end{equation}
Its initial directed azimuth is $0$ when $v_q>v_0$, and $\pi$ when
$v_q<v_0$. In this case, the loxodromic distance is
$\min\{d_q,d_{N}\}$.

If $u_1-u_0$ is not congruent to $0$ modulo $\pi$, no meridian
segment joins $p_0$ and $p_1$. Therefore, the loxodromic distance is
$d_{N}$.

Now suppose that one endpoint is the pole $q$. For every
loxodromic segment of directed azimuth $\theta$ joining
the two endpoints, we have
\begin{equation*}
d=\frac{|v_1-v_0|}{|\cos\theta|},
\end{equation*}
where
$\operatorname{sgn}(\cos\theta) = \operatorname{sgn}(v_1-v_0)$. Consequently, $d\geq|v_1-v_0|$, with equality if and only if the
segment is meridional. Therefore, the loxodromic distance is
$|v_1-v_0|$.

The preceding cases include all nonmeridional segments, all meridian segments between two nonpole points, and all segments having the pole as an endpoint.

\subsubsection*{(d)}
Suppose that $S$ is a sphere-like surface with poles at $v=a$ and $v=b$, and let $L=b-a$. Suppose first that $p_0$ and $p_1$ are distinct nonpole points. Parallel segments and other nonmeridional loxodromic segments are determined as in the cylinder-like case. Let $d_N$ denote the length of the shortest nonmeridional loxodromic segment joining $p_0$ and $p_1$. Thus, $d_N$ is determined by
\eqref{eq:distance_of_loxodrome1} when $v_0=v_1$, and by
\eqref{eq:distance_of_loxodrome2} when $v_0\neq v_1$. In the latter case, the minimum is taken only over
integers $n$ satisfying $u_1-u_0+2\pi n\neq0$.

Meridian segments from $p_0$ to $p_1$ occur when
$u_1-u_0\equiv0\pmod\pi$. We again distinguish two cases.

First, suppose that $u_1-u_0\equiv0\pmod{2\pi}$. Then $p_0$ and $p_1$
lie on the same fixed-longitude meridian branch. The two meridian arcs joining them have lengths
\begin{equation}\label{eq:distance_of_same_longitude_meridians}
d_1=|v_1-v_0|, \quad d_2=2L-|v_1-v_0|.
\end{equation}
All meridional segment lengths are of the form $d_1+2kL$ or $d_2+2kL$, where $k$ is a nonnegative integer. Their initial directed azimuths are $0$ and $\pi$. More precisely, if
$v_1>v_0$, the segment of length $d_1$ has initial directed azimuth
$0$, whereas the segment of length $d_2$ has initial directed azimuth
$\pi$. If $v_1<v_0$, these two directed azimuths are reversed. In this
case, the loxodromic distance is
\begin{equation*}
\min\{d_1,d_2,d_N\}=\min\{d_1,d_N\},
\end{equation*}
since $d_1\leq d_2$.

Next, suppose that $u_1-u_0\equiv\pi\pmod{2\pi}$. Then one meridian
segment passes through the pole at $v=a$, while the other passes
through the pole at $v=b$. Their lengths are
\begin{equation}\label{eq:distance_of_opposite_longitude_meridians}
d_a = |v_0-a|+|v_1-a|, \quad d_b = |b-v_0|+|b-v_1|,
\end{equation}
respectively. All meridional segment lengths are of the form $d_a+2kL$ or $d_b+2kL$, where $k$ is a nonnegative integer. The segment of length $d_a$ has initial directed
azimuth $\pi$, and the segment of length $d_b$ has initial directed
azimuth $0$. In this case, the loxodromic distance is
\begin{equation*}
\min\{d_a,d_b,d_N\}.
\end{equation*}

If $u_1-u_0$ is not congruent to $0$ modulo $\pi$, no meridian segment joins $p_0$ and $p_1$. Therefore, the loxodromic distance is $d_N$.

Finally, suppose that at least one endpoint is a pole. If exactly one endpoint is a pole, the two meridian arcs joining the endpoints have lengths
$$
d_1=|v_1-v_0|,
\qquad d_2=2L-|v_1-v_0|.
$$
All meridional segment lengths are of the form $d_1+2kL$ or $d_2+2kL$, where $k$ is a nonnegative integer.

If $p_0$ is not a pole, the arc of length $d_1$ has initial directed azimuth $0$ when $v_1>v_0$, and $\pi$ when $v_1<v_0$; the arc of length $d_2$ has the opposite initial directed azimuth. If $p_0$ is a pole, both arcs have initial directed azimuth $0$ when $v_0=a$, and $\pi$ when $v_0=b$.

If both endpoints are poles, each fixed-longitude meridian branch gives an arc joining them of length $L$. All meridional segment lengths in this case are $(2k+1)L$, where $k$ is a nonnegative integer. Its initial directed azimuth is $0$ when $v_0=a$, and $\pi$ when $v_0=b$.

In either case, every nonmeridional loxodromic segment of directed azimuth $\theta$ joining the endpoints has length
$$
d=\frac{|v_1-v_0|}{|\cos\theta|},
$$
where
$\operatorname{sgn}(\cos\theta)
=\operatorname{sgn}(v_1-v_0)$.
Thus every such segment has length at least $|v_1-v_0|$. A shortest meridian arc attains this bound, and
additional complete circuits only increase the length. Therefore, the loxodromic distance is $|v_1-v_0|$. Table~\ref{tab:inverse_inverse_summary_refs} summarizes the formulas for the inverse problem.

\section{Geodesics and Loxodromes} \label{sec:geoloxo}
Recall that a $C^2$-curve $\gamma$ on a surface $S$ is called a geodesic if $\gamma''(s)$ is zero or perpendicular to the tangent plane of $S$ at $\gamma(s)$ for every $s$ in its domain. This section investigates the circumstances under which a maximal loxodromic curve on a surface of revolution $S$ is also a geodesic. For our purposes, we first recall some essential preliminary results.

Recall that a meridian, or fixed-longitude meridian branch, is given by $\Phi(u_0,v)$ with $u_0$ fixed, whereas a parallel is given by $\Phi(u,v_0)$ with $v_0$ fixed. Their regular restrictions away from the poles are loxodromes with directed azimuths $0$ or $\pi$ and $\pi/2$ or $3\pi/2$, respectively. These curves should be distinguished from maximal meridians and maximal parallels, which are maximal loxodromic curves parametrized by arc length. When a pole is present, a fixed-longitude meridian branch reaches the pole, whereas a maximal meridian passes smoothly through it along the opposite longitude.

\begin{prop}[Section 9.3 in \cite{Pressley2010}]\label{prop:meridian_parallel} On a surface of revolution with $r>0$,
\begin{enumerate}
    \item Every meridian is a geodesic.
    \item A parallel $\gamma(u) = \Phi(u,v_0)$ is a geodesic if and only if $r'(v_0)= 0$; i.e., $v_0$ is a stationary point of $r$.
\end{enumerate}
\end{prop}
 Proposition~\ref{prop:meridian_parallel} identifies specific families of curves on a surface of revolution that possess geodesic properties. Away from the poles, meridian and parallel branches are special cases of loxodromes. It is therefore natural to investigate when a maximal loxodromic curve is also a geodesic.

\begin{prop}[Clairaut’s Theorem, Section 9.3 in \cite{Pressley2010}]\label{prop:Clairaut_Thm}
    Let $\gamma$ be a unit-speed curve on a surface of revolution $S$, let $\rho : S \to \mathbb{R}$ be the distance of a point of $S$ from the axis of rotation, and let $\theta$ be the angle between $\gamma'$ and the meridians of $S$. If $\gamma$ is a geodesic, then $\rho \sin(\theta)$ is constant along $\gamma$. Conversely, if $\rho \sin(\theta)$ is constant along $\gamma$, and if no part of $\gamma$ is part of some parallel of $S$, then $\gamma$ is a geodesic.
\end{prop}

The distance $\rho$ in Proposition \ref{prop:Clairaut_Thm} is the first component \(r(v)\) of the profile curve in our notation. The next two theorems give conditions for a maximal loxodromic curve to be a geodesic, depending on whether $r$ is constant.

\begin{theorem}\label{thm:loxodrome_on_a_cylinder}
    Every maximal loxodromic curve on a cylinder is a geodesic.
\end{theorem}
\begin{proof}
    Suppose that $S$ is a cylinder, so $r$ is constant. If $\lambda$ is a maximal parallel, then $\lambda$ is a parallel. Since $r'(v)=0$ for every $v$, Proposition~\ref{prop:meridian_parallel} shows that $\lambda$ is a geodesic. Now suppose that $\lambda$ is a maximal loxodromic curve that is not a maximal parallel, and let $\theta$ be the common directed azimuth of its restrictions. Since $r$ and $\sin\theta$ are constant, $r\sin\theta$ is constant along $\lambda$. Moreover, since $\cos\theta\neq0$, no open part of $\lambda$ is contained in a parallel. Therefore, Clairaut's theorem shows that $\lambda$ is a geodesic.

\end{proof}

\begin{theorem}\label{thm:loxodrome_on_a_no-cylinder}
    Let $S$ be a surface of revolution with $r>0$, and assume that $r$ is not constant on any open interval. Then a maximal loxodromic curve is a geodesic if and only if it is either a maximal meridian or a maximal parallel at a stationary point of $r$.
\end{theorem}
\begin{proof}
    By Proposition~\ref{prop:meridian_parallel}, maximal parallels at stationary points of $r$ and maximal meridians are indeed geodesics. We now prove the converse. Let $\lambda(s) = \Phi(u(s),v(s))$ be a geodesic and a maximal loxodromic curve. Since $S$ has no poles, let $\theta$ be the common directed azimuth of the restrictions of $\lambda$. If $\sin\theta=0$, then \eqref{eq:duds} gives $u'(s)=0$, so $\lambda$ is a maximal meridian. If $\cos\theta=0$, then \eqref{eq:dvds} gives $v'(s)=0$, so $\lambda$ is a maximal parallel. Since $\lambda$ is a geodesic, Proposition~\ref{prop:meridian_parallel} implies that $r'(v_0)=0$.

    Suppose for the sake of contradiction that $\sin\theta \neq 0$ and $\cos\theta \neq 0$. By Clairaut's relation, $r(v(s))\sin\theta = C$ for some constant $C$ along $\lambda$. Since $\sin\theta \neq 0$, this implies that $r(v(s))$ is a constant, and we have $r'(v(s))v'(s) = 0$. Since $\cos\theta \neq 0$, we have $v'(s) \neq 0$ for each $s$. Therefore, $r'(v(s)) = 0$ for each $s$. Since $v(s)$ is strictly monotonic, its image on the interior of the parameter interval is a nonempty open interval in the domain of $r$. This implies that $r(v)$ is constant on an open interval, which contradicts the assumption that $r$ is not constant on any open interval. Therefore, $\lambda$ cannot be a geodesic when both $\sin\theta\neq0$ and $\cos\theta\neq0$.
\end{proof}
When we consider maximal loxodromic curves on an arbitrary surface of revolution, only maximal meridians require additional verification at the poles, since they are the only maximal loxodromic curves that can pass through a pole. Every other maximal loxodromic curve either avoids the poles or only reaches a pole as an endpoint. Away from the poles, a maximal meridian is a meridian and hence is a geodesic by Proposition~\ref{prop:meridian_parallel}. It remains to verify the geodesic condition at each pole through which it passes.

\begin{prop}\label{prop:maximal_meridian_geodesic}
Every maximal meridian on a surface of revolution is a geodesic.
\end{prop}
\begin{proof}
Let $\lambda$ be a maximal meridian. If $S$ has no poles, then $\lambda$ is a meridian, and the result follows from Proposition~\ref{prop:meridian_parallel}. Now suppose that $S$ has a pole $q$. Without loss of generality, suppose $v_q=b=0$ and $r(0)=0$. Since $(r'(v))^2+(z'(v))^2=1$ for every $v\in I$, differentiating with respect to $v$ gives
$$
r'(v)r''(v) + z'(v)z''(v) = 0.
$$
Since $|r'(0)| = 1$ and $z'(0) = 0$, we have $r''(0) = 0$. Without loss of generality, suppose that $\lambda(0)=q$. For sufficiently small $\delta>0$, the restriction $\lambda|_{(-\delta,\delta)}$ is given by the following formula. Applying \eqref{eq:meridian_for_1pole}, we have
$$
\lambda(s) =
\begin{cases}
(r(s)\cos u_0, r(s)\sin u_0, z(s)) & \text{if } s \le 0, \\
(-r(-s)\cos u_0, -r(-s)\sin u_0, z(-s)) & \text{if } s > 0.
\end{cases}
$$
The left-hand limit of the difference quotient of $\lambda'$ at $0$ is
\begin{align*}
\lim_{s\to0^-}\frac{\lambda'(s)-\lambda'(0)}{s}
&=(r''(0)\cos u_0,r''(0)\sin u_0,z''(0)),
\end{align*}
whereas the right-hand limit is
\begin{align*}
\lim_{s\to0^+}\frac{\lambda'(s)-\lambda'(0)}{s}
&=(-r''(0)\cos u_0,-r''(0)\sin u_0,z''(0)).
\end{align*}
Since $r''(0)=0$, the two limits coincide. Hence
$$
\lambda''(0)=(0,0,z''(0)).
$$
The tangent plane at the pole is horizontal (parallel to the $xy$-plane). Thus, $\lambda''(0)$ is normal to the tangent plane, so the geodesic condition holds at the pole. Away from the poles, $\lambda$ is a meridian and hence is a geodesic by Proposition~\ref{prop:meridian_parallel}. Therefore, $\lambda$ is a geodesic. Notice that here the meridian is parametrized by arc length, but a different parametrization of the same meridian may not satisfy the definition of a geodesic.
\end{proof}

\section{Mappings Preserving Maximal Loxodromic Curves}
\label{sec:preserving-maximal-loxodromic-curves}

Let $S$ and $\hat S$ be surfaces of revolution. A diffeomorphism $f:S\to\hat S$ is said to \textbf{preserve maximal loxodromic curves} if, for every maximal loxodromic curve $\lambda$ on $S$, an arc-length parametrization of $f\circ\lambda$ is also a maximal loxodromic curve on $\hat S$. If both $f$ and $f^{-1}$ preserve maximal loxodromic curves, then $f$ is called a \textbf{loxodromic equivalence}. In what follows, when we say that $f\circ\lambda$ is a maximal loxodromic curve, we mean that the arc-length parametrization of $f\circ\lambda$ is maximal.

\begin{prop}\label{prop:preserving_loxo_with_no_pole}
Let $f:S\to\hat S$ be a diffeomorphism preserving maximal loxodromic curves. If $\lambda:J\to S$ is a loxodrome and $(f\circ\lambda)(J)$ contains no pole of $\hat S$, then $f\circ\lambda$ is a loxodrome.
\end{prop}
\begin{proof}
Extend an arc-length parametrization of $\lambda$ to a maximal
loxodromic curve. By assumption, its image under $f$ is,
after parametrization by arc length, a maximal loxodromic curve.
Since $(f\circ\lambda)(J)$ contains no pole, the corresponding
restriction is a loxodrome by
Definition~\ref{defi:loxodromic_curve}.
Thus $f\circ\lambda$, being a regular reparametrization
of that restriction, is also a loxodrome.
\end{proof}

We first prove that $f(P)\subseteq \hat{P}$, where $P, \hat{P}$ are the respective sets of poles.

\begin{lemma}\label{lemma:forward_pole_preservation}
Let $f:S\to\hat S$ be a diffeomorphism preserving maximal loxodromic curves. Then $f$ maps every pole of $S$ to a pole of $\hat S$.
\end{lemma}
\begin{proof}
 Let $q$ be a pole of $S$, and choose a nonmeridional maximal loxodromic curve $\lambda$ having $q$ as an endpoint. By assumption, $f\circ\lambda$ is a maximal loxodromic curve on $\hat S$ having $f(q)$ as an endpoint. Suppose that $f(q)$ is not a pole. Let $\mu$ be an arc-length parametrization of $f\circ\lambda$, and let $s_q$ be the endpoint of its domain satisfying $\mu(s_q)=f(q)$. On the side of $s_q$ contained in its domain, the coordinate functions of $\mu$ satisfy \eqref{eq:duds} and \eqref{eq:dvds} for some directed azimuth $\theta$. Since $f(q)$ is not a pole, the coefficients in these equations are smooth near $f(q)$. Hence the solution through $f(q)$ is defined on an open interval containing $s_q$. By uniqueness, this solution agrees with $\mu$ on their common domain near $s_q$. It therefore gives an arc-length-parametrized loxodromic curve on a larger interval that contains $\mu$ as a restriction. This contradicts the maximality of $f\circ\lambda$. Hence $f(q)$ is a pole of $\hat S$.
\end{proof}
Next, we show that $f^{-1}$ maps every loxodrome on $\hat S$
to a loxodrome on $S$.

\begin{lemma}\label{lemma:inverse_on_pole_free_interval}
Let $f:S\to\hat S$ be a diffeomorphism preserving maximal
loxodromic curves. If $\hat\lambda:J\to\hat S$ is a loxodrome
on a nonempty open interval $J$, then
$f^{-1}\circ\hat\lambda$ is a loxodrome.
\end{lemma}

\begin{proof}
Set $\lambda(t)=f^{-1}(\hat\lambda(t))$. Since $f$ is a diffeomorphism, $\lambda$ is smooth and regular. By Lemma~\ref{lemma:forward_pole_preservation}, $\lambda(J)$ contains no pole. Fix $t_0\in J$, and define the arc-length parameter $s(t)=\int_{t_0}^{t}|\lambda'(\tau)|\,d\tau$. Since $s'(t)>0$, let $t=t(s)$ be its inverse. Retaining the notation for the reparametrized curves, write
$$
\lambda(s)=\lambda(t(s)),
\quad
\hat\lambda(s)=\hat\lambda(t(s)),
\quad s\in s(J),
$$
where the right-hand sides use the original parametrizations.
Then $\lambda$ has unit speed, $\hat\lambda$ remains
a loxodrome, and $f\circ\lambda=\hat\lambda$.

Fix $s_0\in s(J)$.
The existence and uniqueness theorem applied to
\eqref{eq:duds} and \eqref{eq:dvds} gives a unit-speed loxodrome $\gamma:K\to S$ satisfying
$$
\gamma(s_0)=\lambda(s_0),
\quad
\gamma'(s_0)=\lambda'(s_0).
$$
Since $f(\gamma(s_0))=\hat\lambda(s_0)$ is not a pole, $s_0$ lies in the inverse image under $f\circ\gamma$ of the set of nonpole points of $\hat S$. The set of nonpole points is open, so its inverse image is open since $f\circ\gamma$ is continuous. Choose an open interval $K_0\subseteq s(J)\cap K$ containing $s_0$ in this inverse image. By Proposition~\ref{prop:preserving_loxo_with_no_pole}, $(f\circ\gamma)|_{K_0}$ is a loxodrome. The chain rule gives
\begin{equation}\label{eq:some_eq_on_inverse_on_pole_free_interval}
    (f\circ\gamma)'(s_0)
=f_{*\lambda(s_0)}\bigl(\lambda'(s_0)\bigr)
=(f\circ\lambda)'(s_0)
=\hat\lambda'(s_0).
\end{equation}
Let $\theta$ be the directed azimuth of $\hat\lambda$. \eqref{eq:some_eq_on_inverse_on_pole_free_interval} implies that $(f\circ\gamma)|_{K_0}$ also has directed azimuth $\theta$. Reparametrize restrictions of $f\circ\gamma$ and $\hat\lambda$ by arc length using increasing changes of parameter, with parameter value $s_0$ at their common point. Their initial unit tangent vectors agree. Since both have directed azimuth $\theta$, their coordinate functions satisfy \eqref{eq:duds} and \eqref{eq:dvds} on $\hat S$ with the same initial conditions. By uniqueness, these arc-length parametrizations agree on an open interval containing $s_0$.

Applying $f^{-1}$, we find that restrictions of $\gamma$ and $\lambda$ agree after an increasing change of parameter fixing $s_0$. Since both have unit speed, this change of parameter has derivative $1$ and is therefore the identity. Hence there exists an open interval $U_{s_0}\subseteq s(J)\cap K_0$ containing $s_0$ such that
$$
\gamma|_{U_{s_0}}=\lambda|_{U_{s_0}}.
$$
Since $s_0$ was arbitrary, every point of $s(J)$ lies in
an open interval on which $\lambda$ has constant directed
azimuth. These azimuths agree whenever the intervals overlap.
Since $s(J)$ is an interval, they agree throughout $s(J)$.
Thus $\lambda$ is a loxodrome.
\end{proof}

Having shown that $f$ maps every pole of $S$ to a pole of $\hat S$, we now prove that every pole of $\hat S$ arises in this way, thereby establishing a bijection between the two sets of poles.
\begin{lemma}\label{lemma:pole_bijection}
Let $f:S\to\hat S$ be a diffeomorphism preserving maximal loxodromic curves. Then $f$ induces a bijection between the poles of $S$ and the poles of $\hat S$.
\end{lemma}

\begin{proof}
By Lemma~\ref{lemma:forward_pole_preservation}, $f$ maps every pole of $S$ to a pole of $\hat S$. It remains to show that $f^{-1}(\hat q)$ is a pole of $S$ for every pole $\hat q$ of $\hat S$.

Choose a nonmeridional maximal loxodromic curve $\hat\lambda:J\to\hat S$ having $\hat q$ as an endpoint. By Corollary~\ref{cor:topo_of_loxodrome_on_PS} or \ref{cor:topo_of_loxodrome_on_SS}, $\hat\lambda|_{\operatorname{int}J}$ is a loxodrome. Lemma~\ref{lemma:inverse_on_pole_free_interval} therefore implies that $(f^{-1}\circ\hat\lambda)|_{\operatorname{int}J}$ is a loxodrome.

Suppose that $f^{-1}(\hat q)$ is not a pole.
Let $\mu:K\to S$ be an arc-length parametrization of
$f^{-1}\circ\hat\lambda$.
The curve $\hat\lambda$ reaches $\hat q$ in finite length.
Since the differential of $f^{-1}$ is bounded on a
neighborhood of $\hat q$, the corresponding portion
of $f^{-1}\circ\hat\lambda$ also has finite length.
Thus $K$ has a finite endpoint $s_{\hat q}$ satisfying
$$
\mu(s_{\hat q})=f^{-1}(\hat q).
$$

On $\operatorname{int}K$, the coordinate functions of $\mu$
satisfy \eqref{eq:duds} and \eqref{eq:dvds}
for some directed azimuth $\theta$.
Since $f^{-1}(\hat q)$ is not a pole, $r$ is positive
at its profile parameter.
The right-hand sides of these equations are therefore
smooth on an open coordinate neighborhood of
$f^{-1}(\hat q)$.
The existence and uniqueness theorem extends $\mu$
through $s_{\hat q}$.
Choose a maximal loxodromic curve $\gamma$
containing this extension as a restriction.

By assumption, an arc-length parametrization of
$f\circ\gamma$ is a maximal loxodromic curve on $\hat S$.
Choose its direction and parameter origin so that it
agrees with $\hat\lambda$ on $J$.
Since $\gamma$ extends $\mu$ through $s_{\hat q}$,
this parametrization extends $\hat\lambda$ beyond
its endpoint, contradicting maximality.
Therefore, $f^{-1}(\hat q)$ is a pole of $S$.

Together with Lemma~\ref{lemma:forward_pole_preservation}
and the injectivity of $f$, this proves the assertion.
\end{proof}

By Lemma~\ref{lemma:pole_bijection} and Proposition~\ref{prop:preserving_loxo_with_no_pole}, we conclude that $f$ maps every loxodrome to a loxodrome. Using the preceding lemmas, we show that every diffeomorphism preserving maximal loxodromic curves is a loxodromic equivalence. 

\begin{prop}\label{prop:loxodromic-equivalence}
Let $f:S\to\hat S$ be a diffeomorphism preserving maximal loxodromic curves. Then $f^{-1}$ also preserves maximal loxodromic curves. Hence $f$ is a loxodromic equivalence.
\end{prop}
\begin{proof}
Let $\hat\lambda:J\to\hat S$ be a maximal loxodromic curve, and set
$\lambda=f^{-1}\circ\hat\lambda$. Since $f^{-1}$ is a diffeomorphism, $\lambda$ is smooth and regular on $\operatorname{int}J$.

Let $J_0\subseteq\operatorname{int}J$ be a nonempty open interval such that $\lambda(J_0)$ contains no pole of $S$. By Lemma~\ref{lemma:pole_bijection}, $\hat\lambda(J_0)$ contains no pole of $\hat S$. Lemma~\ref{lemma:inverse_on_pole_free_interval} therefore shows that $\lambda|_{J_0}$ is a loxodrome. Hence $\lambda$ is a loxodromic curve.

Let $\eta:L\to S$ be an arc-length parametrization of $\lambda$. Suppose that $\eta$ is not maximal. Then there are an interval $K$ properly containing $L$ and an arc-length-parametrized loxodromic curve $\gamma:K\to S$ such that $\gamma|_L=\eta$. 

Choose a maximal loxodromic curve $\nu$ such that $\nu|_K=\gamma.$
Since $f$ preserves maximal loxodromic curves and $\nu$
is maximal, an arc-length parametrization of $f\circ\nu$
is a maximal loxodromic curve on $\hat S$. Since $f\circ\eta$ is a regular reparametrization of $\hat\lambda$, we may choose this arc-length parametrization so that it agrees with $\hat\lambda$ on $J$. Since $f$ is a diffeomorphism, $f\circ\gamma$ is regular, so passage to arc length preserves the proper inclusion $L\subseteq K$. Consequently, the domain of the chosen arc-length parametrization properly contains $J$. Thus, $\hat\lambda$ is the restriction of an arc-length-parametrized loxodromic curve defined on a larger interval. This contradicts the maximality of $\hat\lambda$.

Therefore, $\eta$ is maximal. Since $\hat\lambda$ was arbitrary, $f^{-1}$ preserves maximal loxodromic curves. Hence $f$ is a loxodromic equivalence.
\end{proof}

We next show that such diffeomorphisms preserve the classification of surfaces of revolution.

\begin{theorem}\label{prop:classification_invariance'}
    Suppose $f:S\to\hat S$ is a diffeomorphism preserving maximal loxodromic curves. If $S$ is a cylinder-like surface, then $\hat S$ is also cylinder-like. Similar correspondences hold for torus-like, paraboloid-like, and sphere-like surfaces.
\end{theorem}

\begin{proof}
    According to Lemma~\ref{lemma:pole_bijection}, the surfaces $S$ and $\hat S$ have the same number of poles. Therefore, if $S$ is paraboloid-like or sphere-like, then $\hat S$ is paraboloid-like or sphere-like, respectively. It remains to distinguish the cylinder-like and torus-like cases. Suppose that $S$ is torus-like. Then there exists a dense maximal loxodromic curve $\lambda$ on $S$. Since $f$ is a homeomorphism, $f\circ\lambda$ is dense in $\hat S$. Since $f\circ\lambda$ is also a maximal loxodromic curve, $\hat S$ is torus-like. Finally, suppose that $S$ is cylinder-like. Since $\hat S$ has no poles, it is either cylinder-like or torus-like. If $\hat S$ were torus-like, Proposition~\ref{prop:loxodromic-equivalence} would imply that $f^{-1}:\hat S\to S$ preserves maximal loxodromic curves. Applying the preceding torus-like result to $f^{-1}$ would then imply that $S$ is torus-like, a contradiction. Therefore, $\hat S$ is cylinder-like.
\end{proof}

In the remainder of the article, we examine whether a loxodromic equivalence $f$ maps parallels and meridians of a given surface $S$ to parallels and meridians of $\hat S$, depending on the type of $S$ and hence of $\hat S$ by Theorem ~\ref{prop:classification_invariance'}. The following four subsections, also summarized in Table~\ref{tab:preservation_properties}, treat the four types separately. For cylinder-like surfaces, $f$ preserves parallels but does not necessarily preserve meridians. For paraboloid-like and sphere-like surfaces, $f$ preserves both parallels and meridians. For torus-like surfaces, $f$ does not necessarily preserve either parallels or meridians.

In general, finding a precise differential formula for all loxodromic equivalences is difficult. To obtain a more useful result, we introduce the following stronger condition, which will also apply to the surfaces considered in the following subsections.

Let $P$ and $\hat P$ denote the sets of poles of $S$ and $\hat S$, respectively, and set $S^*=S\setminus P$ and $\hat S^*=\hat S\setminus\hat P$. By Lemma~\ref{lemma:pole_bijection}, a loxodromic equivalence maps $P$ bijectively onto $\hat P$, and hence $f(S^*)=\hat S^*$. A loxodromic equivalence $f:S\to\hat S$ is called \textbf{azimuth-consistent} if, for every $\theta\in[0,2\pi)$, there exists a unique $\hat\theta\in[0,2\pi)$, depending only on $\theta$, such that the unit vector field $f_*W_\theta/|f_*W_\theta|$ has constant directed azimuth $\hat\theta$ on $\hat S^*$. That is, for every $p\in S^*$,
$$
\frac{f_{*p}(W_\theta|_p)}
     {|f_{*p}(W_\theta|_p)|}
=
\hat W_{\hat\theta}|_{f(p)},
$$
where $W_\theta=\sin\theta\,e_u+\cos\theta\,e_v$ on $S^*$ and $\hat W_{\hat\theta} =\sin\hat\theta\,\hat e_u+\cos\hat\theta\,\hat e_v$ on $\hat S^*$ are defined by \eqref{eq:Vector_field_loxo_arc}. 
In particular, $\hat\theta$ depends only on $\theta$
and is independent of $p\in S^*$.

\subsection{Cylinder-Like Surfaces}
Let $S$ and $\hat S$ be cylinder-like surfaces, and suppose that $f:S\to\hat S$ is a loxodromic equivalence. Let $\lambda$ be a maximal parallel at $v=v_0$ on $S$. Since $\lambda$ is closed, $f\circ\lambda$ is also closed. Moreover, since $f$ is a loxodromic equivalence, $f\circ\lambda$ is a maximal loxodromic curve on $\hat S$. By Proposition~\ref{prop:topo_of_loxodrome_on_CS}, it must be a maximal parallel. Therefore, we have the following result.
\begin{prop}\label{prop:map_to_para}
    Let $S$ and $\hat S$ be cylinder-like surfaces. Suppose that $f:S\to\hat S$ is a loxodromic equivalence. If $\lambda$ is a maximal parallel at $v=v_0$ on $S$, then $f\circ\lambda$ is a maximal parallel on $\hat S$.
\end{prop}
  Using equations (\ref{eq:naturalTV_u}) and (\ref{eq:naturalTV_v}), we define the orthonormal frames $\{e_u, e_v\}$ and $\{\hat{e}_u, \hat{e}_v\}$ for the tangent bundles of $S$ and $\hat{S}$, respectively, by
$$
e_u = \frac{1}{r}\partial_u, \quad e_v = \partial_v, \qquad \hat{e}_u = \frac{1}{\hat{r}}\partial_{\hat{u}}, \quad \hat{e}_v = \partial_{\hat{v}}.
$$
By the definition of an azimuth-consistent diffeomorphism, $f$ is azimuth-consistent on a cylinder-like surface if and only if, for each unit vector field $X$ making a constant angle with $e_v$, the unit vector field $f_*X/|f_*X|$ also makes a constant angle with $\hat{e}_v$.

According to Proposition~\ref{prop:map_to_para}, for each $p\in S$ and $q = f(p)$, we have
\begin{align}\label{eq:eu_to_hateu}
     f_*(e_u\big\rvert_{p})= k(p) \hat{e}_u\big\rvert_{q}
\end{align}
    for some nonvanishing function $k \in C^\infty(S)$.

\begin{theorem}\label{thm:LP_and_OGmatrax}
    Let $S$ and $\hat S$ be cylinder-like surfaces. A diffeomorphism $f:S\to\hat S$ is azimuth-consistent if and only if its differential satisfies
    \begin{align}
    f_*(e_u\big\rvert_{p})&=   k(p) \hat{e}_u\big\rvert_{q},\label{formula:differential_eu}\\
    f_*(e_v\big\rvert_{p})&=  c_1 k(p) \hat{e}_u\big\rvert_{q}+c_2k(p) \hat{e}_v\big\rvert_{q},\label{formula:differential_ev}
    \end{align}
    where $p\in S$, $q = f(p)$, $k$ is a smooth nonvanishing function on $S$, and $c_1, c_2$ are constants with $c_2\neq 0$.
\end{theorem}
\begin{proof}
    Suppose $f$ is azimuth-consistent. Applying (\ref{eq:eu_to_hateu}), we have $f_*(e_u\rvert_{p})= k(p)\hat{e}_u\rvert_{q}$. Let $f_*(e_v\rvert_{p})= k_1(p) \hat{e}_u\rvert_{q}+k_2(p) \hat{e}_v\rvert_{q}$, where $k_1, k_2 \in C^\infty(S)$. Since $f$ is a diffeomorphism, the functions $k$ and $k_2$ are nonvanishing. Let $X =\sin(\theta){e}_u + \cos(\theta){e}_v $ be a unit vector field of constant directed azimuth $\theta$, where $\cos\theta\neq0$. We have
    \begin{equation}\label{eq:pushforward_X}
        f_*(X\big\rvert_{p}) = (k_1(p)\cos(\theta) +k(p)\sin(\theta))\hat{e}_u\rvert_{q}+k_2(p)\cos(\theta) \hat{e}_v\rvert_{q}.
    \end{equation}
    From (\ref{eq:pushforward_X}), the tangent of the directed azimuth $\hat{\theta}$ of $f_*(X)$ is given by
    \begin{equation}\label{eq:azimuth_relation}
    \tan(\hat\theta) = \frac{k_1(p)\cos(\theta) +k(p)\sin(\theta)}{k_2(p)\cos(\theta)} = \frac{k_1(p)}{k_2(p)}+\frac{k(p)}{k_2(p)}\tan(\theta)
    \end{equation}

    for each $p \in S$. Since $f$ is azimuth-consistent, the directed azimuth $\hat{\theta}$ is constant and independent of $p$ for each fixed directed azimuth $\theta$. Evaluating (\ref{eq:azimuth_relation}) at $\theta = 0$ and $\theta = \pi/4$ yields
    \begin{align*}
        \frac{k_1(p)}{k_2(p)} &= \tan(\hat\theta(0)),\\
        \frac{k(p)}{k_2(p)} &= \tan(\hat\theta(\pi/4)) - \tan(\hat\theta(0)).
    \end{align*}
    Therefore, $k_1/k_2$ and $k/k_2$ are constants. Since $k$ and $k_2$ are nonzero, $k_1/k$ and $k_2/k$ are also constants. Hence there exist constants $c_1$ and $c_2$, with $c_2\neq0$, such that (\ref{formula:differential_eu}) and (\ref{formula:differential_ev}) hold.

    Conversely, suppose the differential of $f$ takes the given forms (\ref{formula:differential_eu}) and (\ref{formula:differential_ev}). Let $X =\sin(\theta){e}_u + \cos(\theta){e}_v$ be a unit vector field of constant directed azimuth $\theta$. We have
    $$
    f_*(X\big\rvert_{p}) = \big(c_1 k(p)\cos(\theta) +  k(p)\sin(\theta)\big)\hat{e}_u\big\rvert_{q}+c_2k(p)\cos(\theta) \hat{e}_v\big\rvert_{q}.
    $$
    If $\cos(\theta)\neq 0$, then the tangent of the directed azimuth $\hat{\theta}$ of $f_*(X)$ is given by
    \begin{align*}
    \tan(\hat{\theta}) = \frac{ c_1k(p)\cos(\theta) +  k(p)\sin(\theta)}{c_2k(p)\cos(\theta)} = \frac{c_1}{c_2} +\frac{1}{c_2}\tan(\theta).
    \end{align*}
    Since $c_1, c_2$, and $\theta$ are constants, $\hat{\theta}$ is constant. If $\cos(\theta) = 0$, then
    $$
    f_*(X\big\rvert_{p}) =  \pm k(p)\hat{e}_u\rvert_{q}.
    $$
    This implies $\cos \hat{\theta}=0$, meaning that the directed azimuth $\hat{\theta}$ is also constant. Since $S$ is connected and $k$ is nonvanishing, the sign of $k$ is constant. Therefore, in both cases, the directed azimuth $\hat\theta$, and not merely $\tan\hat\theta$, is constant. Thus, for every constant directed azimuth $\theta$, the vector field $f_*X$ has a constant directed azimuth $\hat\theta$ depending only on $\theta$. Therefore, $f$ maps every loxodrome to a loxodrome. Moreover, the differential of $f^{-1}$ has the same form. Indeed, relative to the orthonormal frames, the matrix representation of the differential $f_*$ at $p \in S$ is given by $
k(p)
\begin{pmatrix}
1&c_1\\
0&c_2
\end{pmatrix},
$
and hence the matrix representation of the differential $(f^{-1})_{*q}$ at $q=f(p)\in\hat S$ is given by
$
{1}/{k(p)}
\begin{pmatrix}
1&-{c_1}/c_2\\
0&{1}/{c_2}
\end{pmatrix}.
$
Therefore, both $f$ and $f^{-1}$ map loxodromes to loxodromes. Let $\lambda:J\to S$ be a maximal loxodromic curve, and let $\eta:L\to\hat S$ be an arc-length parametrization of $f\circ\lambda$. Then $\eta$ is a loxodromic curve. Suppose that $\eta$ is not maximal. Then there are an interval $K$ containing $L$ with $K\neq L$ and an arc-length-parametrized loxodromic curve $\gamma:K\to\hat S$ such that $\gamma|_L=\eta$. The curve $f^{-1}\circ\gamma$ is regular and loxodromic. Choose its arc-length parametrization so that it agrees with $\lambda$ on $J$. 
The change to arc length is strictly increasing, so the proper inclusion $L\subset K$ remains proper after reparametrization. Thus the domain of this arc-length parametrization properly contains $J$, contradicting the maximality of $\lambda$. Hence $f$ preserves maximal loxodromic curves. Proposition~\ref{prop:loxodromic-equivalence} shows that $f$ is a loxodromic equivalence. Together with the preceding relation between the directed azimuths, this proves that $f$ is azimuth-consistent.
\end{proof}

\begin{rmk}\label{rmk:LP_and_Nmatrax}
    Let $(u,v)$ and $(\hat{u}, \hat{v})$ be the coordinates of $p$ and $q=f(p)$, respectively. Using (\ref{formula:differential_eu}) and (\ref{formula:differential_ev}), we can express the differential of \(f\) in the coordinate bases as
    \begin{align}
        f_*\big(\partial_{u}\big|_{p}\big) &=  k(p) \frac{r(v)}{\hat{r}(\hat{v})} \partial_{\hat{u}}\big|_{q}, \\[1.5ex]
        f_*\big(\partial_{v}\big|_{p}\big) &= c_1 \frac{k(p)}{\hat{r}(\hat{v})} \partial_{\hat{u}}\big|_{q} +c_2 k(p) \partial_{\hat{v}}\big|_{q}.
    \end{align}
    This expression explicitly relates the differential to the radial functions $r$ and $\hat{r}$. Furthermore, its matrix representation is given by
    \begin{align}\label{eq:matrix_form}
        J=\begin{pmatrix}
        \frac{\partial \hat{u}}{\partial u} & \frac{\partial \hat{u}}{\partial v}\\
        \frac{\partial \hat{v}}{\partial u} & \frac{\partial \hat{v}}{\partial v}
        \end{pmatrix}
        =
        \frac{k(p)}{\hat{r}(\hat{v})}
\begin{pmatrix}
     r(v) & c_1 \\
    0 & c_2\hat{r}(\hat{v})
\end{pmatrix}.
    \end{align}

\end{rmk}

\begin{rmk}\label{rmk:azimuth_map}
    For a maximal loxodromic curve $\lambda$ with azimuth $\theta\not\equiv \pi/2\pmod\pi$, the azimuth $\hat{\theta}$ of $f\circ\lambda$ satisfies
    \begin{align}\label{eq:tan_of_map}
        \tan(\hat{\theta}) = \frac{c_1}{c_2} +\frac{1}{c_2} \tan(\theta).
    \end{align}
\end{rmk}

Moreover, according to (\ref{eq:tan_of_map}), if $\lambda$ is a meridian on $S$, its image $f\circ \lambda$ may not be a meridian on $\hat{S}$ since the constant $c_1$ may be nonzero.

\begin{eg}\label{eg:cylinder_to_cylinder}
    Let $S$ and $\hat{S}$ be unit cylinders in $\mathbb{R}^3$, parametrized by $\Phi(u,v) = (\cos u, \sin u, v)$. Define a diffeomorphism $f: S \to \hat{S}$ by $f(x,y,z) = (x \cos(cz) - y \sin(cz), x \sin(cz) + y \cos(cz), z)$ for some constant $c \neq 0$. This map corresponds to the coordinate transformation
    \begin{align*}
        \hat{u}(u,v) = u + cv,\quad \hat{v}(u,v) = v .
    \end{align*}
    The differential $f_*$ is represented by the Jacobian matrix
    $$
    J = \begin{pmatrix}
        \frac{\partial \hat{u}}{\partial u} & \frac{\partial \hat{u}}{\partial v}\\
        \frac{\partial \hat{v}}{\partial u} & \frac{\partial \hat{v}}{\partial v}
    \end{pmatrix}
    =
    \begin{pmatrix}
        1 & c\\
        0 & 1
    \end{pmatrix}.
    $$
    Comparing this with (\ref{eq:matrix_form}), we find $k \equiv 1$ on $S$, $c_1 = c$, and $c_2 = 1$. Thus, the mapping $f$ is azimuth-consistent. For any meridian $\lambda$ on $S$ with azimuth $\theta=0$ modulo $\pi$, (\ref{eq:tan_of_map}) gives $\tan \hat{\theta} = c$, where $\hat{\theta}$ is the azimuth of the image curve $f\circ\lambda$. Since $c$ is nonzero, $\hat{\theta}$ is not congruent to $0$ modulo $\pi$; therefore, the image curve $f\circ\lambda$ is a loxodrome but not a meridian.
\end{eg}

\begin{cor}\label{cor:azimuth_meridian_preservation}
    An azimuth-consistent diffeomorphism $f: S \to \hat{S}$ preserves meridians if and only if $c_1 = 0$.
\end{cor}

\begin{proof}
    Let $\lambda$ be a meridian on $S$. By Remark \ref{rmk:azimuth_map}, the azimuth $\hat{\theta}$ of the image curve $f \circ \lambda$ satisfies $\tan \hat{\theta} = c_1/c_2$. The image curve $f \circ \lambda$ is a meridian on $\hat{S}$ if and only if $\tan \hat{\theta} = 0$, which is equivalent to $c_1 = 0$.
\end{proof}
 An azimuth-consistent diffeomorphism satisfying $c_1 = 0$ preserves both meridians and parallels. We consider whether a mapping $f:S \to \hat{S}$ between cylinder-like surfaces that preserves both meridians and parallels is necessarily azimuth-consistent. In general, this does not hold. Consider the differential $f_*$. We have
\begin{align*}
    f_*(e_u\big\rvert_{p})&= k_1(p) \hat{e}_u\big\rvert_{q},\\
    f_*(e_v\big\rvert_{p})&=  k_2(p) \hat{e}_v\big\rvert_{q},
\end{align*}
where $p\in S$, $q = f(p)$, and $k_1$, $k_2$ are smooth nonvanishing functions on $S$. By Theorem \ref{thm:LP_and_OGmatrax}, such a map is azimuth-consistent if and only if the ratio $k_2/k_1$ is constant on $S$. Therefore, preservation of meridians and parallels alone does not imply azimuth-consistency. The following example shows that such a map may even fail to preserve maximal loxodromic curves.
\begin{eg}\label{eg:meridian_parallel_not_azimuth_consistent}
    The following example shows that a diffeomorphism preserving meridians and parallels need not be azimuth-consistent. Let $S = \hat{S} = \{(\cos u,\sin u,v):u\in\mathbb{R}, v>0\}$ be the upper unit cylinder in $\mathbb{R}^3$, parametrized by $\Phi(u,v) = (\cos u, \sin u, v)$. Define a diffeomorphism $f: S \to \hat{S}$ by $f(x,y,z) = (x, y, z^3)$, which corresponds to the coordinate transformation $\hat{u}(u,v) = u$ and $\hat{v}(u,v) = v^3$. Observe that $f$ preserves both meridians and parallels. The differential $f_*$ is represented by the Jacobian matrix
    $$
    J = \begin{pmatrix}
        \frac{\partial \hat{u}}{\partial u} & \frac{\partial \hat{u}}{\partial v}\\
        \frac{\partial \hat{v}}{\partial u} & \frac{\partial \hat{v}}{\partial v}
    \end{pmatrix}
    =
    \begin{pmatrix}
        1 & 0\\
        0 & 3v^2
    \end{pmatrix}.
    $$
However, this matrix cannot have the form \eqref{eq:matrix_form} with constants $c_1$ and $c_2$, because the ratio of its diagonal entries is not constant. Therefore, by Theorem \ref{thm:LP_and_OGmatrax}, $f$ is not azimuth-consistent. In fact, $f$ does not preserve maximal loxodromic curves. Indeed, consider the loxodrome $\lambda(v)=\Phi(v,v)$ for $v>0$, whose arc-length parametrization is maximal. Its image is $(f\circ\lambda)(v)=\hat\Phi(v,v^3)$. Writing $\hat v=v^3$, we obtain $\hat u(\hat v)=\hat v^{1/3}$ and hence $\frac{d \hat u}{d \hat v}=\frac{1}{3\hat v^{2/3}}$, which is not constant. Thus, by Corollary~\ref{cor:uv}, $f\circ\lambda$ does not have constant azimuth and is not a loxodrome.
\end{eg}

A local diffeomorphism $f: S \to \hat{S}$ is called \textbf{conformal} if it preserves the angle of intersection between any two intersecting curves on $S$. Equivalently, in terms of local coordinate patches, $f$ is conformal if and only if the first fundamental forms of $S$ and $\hat{S}$ are proportional (see \cite{Pressley2010}).
\begin{prop}\label{prop:LP_and_Cfm}
    An azimuth-consistent diffeomorphism $f:S\to \hat{S}$ is conformal if and only if $c_1=0$ and $c_2 = \pm 1$. Moreover, if $f$ is conformal, the mapping satisfies
    $$
    \langle f_*X_p, f_*Y_p \rangle_{\hat{S}} = k^2(p)\langle X_p,Y_p \rangle_S
    $$
    for each $p \in S$ and $X_p, Y_p \in T_p S$.
\end{prop}

\begin{proof}
    Let $f: S \to \hat{S}$ be an azimuth-consistent diffeomorphism. Let $E, F, G$ be the coefficients of the first fundamental form of $S$ with respect to the parametrization $\Phi$, and let $\hat{E}, \hat{F}, \hat{G}$ be the coefficients of the first fundamental form of $\hat{S}$ with respect to the parametrization $f\circ\Phi$. We have
    \begin{align*}
        E(u,v) = r^2(v), \quad F(u,v) = 0, \quad G(u,v) = 1,
    \end{align*}
    and
    \begin{align*}
        \hat{E}(u,v) &= \langle f_*\partial_u, f_*\partial_u \rangle_{\hat{S}} =  k^2(p) r^2(v), \\
        \hat{F}(u,v) &= \langle f_*\partial_u, f_*\partial_v \rangle_{\hat{S}} = c_1 k^2(p) r(v), \\
        \hat{G}(u,v) &= \langle f_*\partial_v, f_*\partial_v \rangle_{\hat{S}} = (c_1^2+c_2^2)k^2(p).
    \end{align*}
    Suppose $f$ is conformal. Then there exists a positive smooth function $\mu$ on $S$ such that $\hat{E} = \mu E$, $\hat{F} = \mu F$, and $\hat{G} = \mu G$. Since $\hat{E} = \mu E$, we conclude that $k^2 = \mu$. Since $F = 0$, we have $\hat{F} = 0$ so that
    $c_1 k^2(p) r(v) = 0$. Since $k^2$ is a positive smooth function and $r(v) > 0$, it follows that $c_1 = 0$. Substituting $c_1 = 0$ into the equation $\hat{G} = \mu G$, we get $c^2_2k^2(p) = \mu(p)$. Therefore, $c_1 = 0$ and $c_2 = \pm 1$.

    Conversely, suppose $c_1 = 0$ and $c_2 = \pm 1 $. Substituting these constants into the equations for $\hat{E}, \hat{F}, \hat{G}$, we get $\hat{E} = k^2(p) E$, $\hat{F} = 0 = k^2(p) F$, and $\hat{G} = k^2(p) G$. Since the first fundamental forms are proportional with factor $k^2(p)$, the mapping $f$ is conformal.
\end{proof}
\begin{eg}
    (Mercator Projection.) Let $S$ be the unit sphere minus its two poles in $\mathbb{R}^3$, parametrized by $\Phi(u,v) = (\cos v \cos u, \cos v \sin u, \sin v)$, where $v \in (-\pi/2, \pi/2)$. Let $\hat{S}$ be the unit cylinder parametrized by $\hat{\Phi}(\hat{u}, \hat{v}) = (\cos \hat{u}, \sin \hat{u}, \hat{v})$. The Mercator projection is a diffeomorphism $f: S \to \hat{S}$ given by
    $$
        f(x, y, z) = \left( \frac{x}{\sqrt{x^2+y^2}}, \frac{y}{\sqrt{x^2+y^2}}, \frac{1}{2} \ln \left( \frac{1+z}{1-z} \right) \right),
    $$
    which corresponds to the coordinate transformation $\hat{u}(u,v) = u $ and $\hat{v}(u,v) = \ln(\sec(v)+\tan(v))$. The differential $f_*$ is represented by the Jacobian matrix
    $$
    J = \begin{pmatrix}
            \frac{\partial \hat{u}}{\partial u} & \frac{\partial \hat{u}}{\partial v}\\
            \frac{\partial \hat{v}}{\partial u} & \frac{\partial \hat{v}}{\partial v}
        \end{pmatrix}
        =
        \begin{pmatrix}
            1 & 0\\
            0 & \frac{1}{\cos v}
        \end{pmatrix}.
    $$
    Comparing this with (\ref{eq:matrix_form}), we identify $k(p) = \frac{1}{\cos v}$, $c_1 = 0$, and $c_2 = 1$. This implies that the mapping $f$ is azimuth-consistent by Theorem \ref{thm:LP_and_OGmatrax}. Consequently, it is conformal by Proposition~\ref{prop:LP_and_Cfm}.
\end{eg}
An azimuth-consistent diffeomorphism need not be conformal. For example, the diffeomorphism in Example~\ref{eg:cylinder_to_cylinder} has $c_1\neq0$ and is therefore not conformal by Proposition~\ref{prop:LP_and_Cfm}.

A smooth map $f: S \to \hat{S}$ is called a local isometry if it maps any curve on $S$ to a curve of the same length on $\hat{S}$. For a local diffeomorphism, this property is equivalent to the requirement that the first fundamental forms of $S$ and $\hat{S}$ are identical (see \cite{Pressley2010}). Since a conformal mapping has proportional first fundamental forms with a factor of $k^2(p)$, it becomes an isometry if and only if this factor is exactly $1$.

\begin{cor}\label{cor:azimuth_isometry}
    An azimuth-consistent diffeomorphism $f:S\to \hat{S}$ is an isometry if and only if $f$ is conformal with $k^2 \equiv 1$ on $S$.
\end{cor}

\subsection{Paraboloid-Like Surfaces}
By an argument analogous to that of Proposition~\ref{prop:map_to_para}, a loxodromic equivalence $f:S\to\hat S$ between paraboloid-like surfaces maps maximal parallels to maximal parallels. By Lemma~\ref{lemma:pole_bijection}, the pole $q$ of $S$ is mapped to the pole $\hat q=f(q)$ of $\hat S$. The maximal meridians of $S$ are the only maximal loxodromic curves that pass through $q$. Since $f$ is a loxodromic equivalence, the image of each maximal meridian is a maximal loxodromic curve that passes through $\hat q$. It is therefore a maximal meridian on $\hat S$.

\begin{cor}\label{cor:map_to_parallel_and_meridian_for_Paraboloid}
A loxodromic equivalence $f:S\to\hat S$ between paraboloid-like surfaces maps maximal parallels to maximal parallels and maximal meridians to maximal meridians.
\end{cor}

\begin{theorem}\label{thm:azimuth_consistent_paraboloid}
Let $S$ and $\hat S$ be paraboloid-like surfaces with poles $q$ and
$\hat q$, respectively. A diffeomorphism $f:S\to\hat S$ satisfying
$f(q)=\hat q$ is azimuth-consistent if and only if its differential on
$S^*=S\setminus\{q\}$ satisfies
\begin{align}
f_*\left(e_u\big|_p\right)
&=k(p)\hat e_u\big|_{f(p)},\label{eq:push_u_Paraboloid}\\
f_*\left(e_v\big|_p\right)
&=c k(p)\hat e_v\big|_{f(p)},\label{eq:push_v_Paraboloid}
\end{align}
where $p\in S^*$, $k$ is a smooth nonvanishing function on $S^*$,
and $c$ is a nonzero constant.
\end{theorem}

\begin{proof}
Suppose first that $f$ is azimuth-consistent. By Corollary~\ref{cor:map_to_parallel_and_meridian_for_Paraboloid}, $f$ maps maximal meridians to maximal meridians and maximal parallels to maximal parallels. Therefore, by the same argument as in Theorem~\ref{thm:LP_and_OGmatrax}, its differential on $S^*$ has the stated form.

Conversely, suppose that $f$ satisfies \eqref{eq:push_u_Paraboloid} and \eqref{eq:push_v_Paraboloid}. Then both $f$ and $f^{-1}$ map loxodromes to loxodromes away from the poles, and the directed azimuth of the image depends only on the directed azimuth of the original loxodrome. The two equations also show that both maps preserve meridional directions. Since $f$ and $f^{-1}$ are smooth at the poles, they map maximal meridians passing through a pole to maximal meridians passing through the corresponding pole. Every other nonparallel maximal loxodromic curve that reaches a pole has that pole as an endpoint. Since $f(q)=\hat q$, its image also has the corresponding pole as an endpoint. It follows that both $f$ and $f^{-1}$ map loxodromic curves to loxodromic curves.

Let $\lambda:J\to S$ be a maximal loxodromic curve, and let $\eta:L\to\hat S$ be an arc-length parametrization of $f\circ\lambda$. Suppose that $\eta$ is not maximal. Then there are an interval $K$ containing $L$ with $K\neq L$ and an arc-length-parametrized loxodromic curve $\gamma:K\to\hat S$ such that $\gamma|_L=\eta$. The curve $f^{-1}\circ\gamma$ is regular and loxodromic. Choose its arc-length parametrization so that it agrees with $\lambda$ on $J$. The change to arc length is strictly increasing, so the proper inclusion $L\subset K$ remains proper after reparametrization. Thus the domain of this arc-length parametrization properly contains $J$, contradicting the maximality of $\lambda$. Hence $f$ preserves maximal loxodromic curves. Proposition~\ref{prop:loxodromic-equivalence} shows that $f$ is a loxodromic equivalence. Together with the preceding relation between the directed azimuths, this shows that $f$ is azimuth-consistent.
\end{proof}

For an azimuth-consistent diffeomorphism between paraboloid-like surfaces, smoothness at the pole further forces $|c|=1$. Indeed, let $X_q$ be an arbitrary unit tangent vector at $q$. Approach $q$ along a meridian whose tangent at $q$ is $X_q$. Up to signs, the vectors $e_v$ and $e_u$ approach $X_q$ and $RX_q$, respectively, where $R$ denotes rotation by $\pi/2$ in the tangent plane. Taking the limit in $\frac{|f_\ast(e_v)|}{|f_\ast(e_u)|}=|c|$ gives
$$
\frac{|f_{\ast q}(X_q)|}{|f_{\ast q}(RX_q)|}=|c|.
$$
Now choose another meridian whose tangent at $q$ is $RX_q$. Along this meridian, up to signs, $e_v$ approaches $RX_q$, while $e_u$ approaches $R(RX_q)=-X_q$. Applying the same argument yields
$$
\frac{|f_{\ast q}(RX_q)|}{|f_{\ast q}(X_q)|}=|c|.
$$
Multiplying the preceding two equalities gives $|c|^2=1$, and hence $c=\pm1$. Therefore, $f$ is conformal on $S\setminus\{q\}$. Since $f$ is smooth and conformality of its differential is preserved under taking limits, $f$ is also conformal at $q$. Consequently, every azimuth-consistent diffeomorphism between paraboloid-like surfaces is conformal.

\subsection{Sphere-Like Surfaces}
Let $S$ and $\hat S$ be sphere-like surfaces, and let $f:S\to\hat S$ be a loxodromic equivalence. Since $f$ induces a bijection between the poles of $S$ and those of $\hat S$, the image of every maximal meridian passes through the two poles of $\hat S$ and is therefore a maximal meridian.

Let $\lambda$ be a maximal parallel on $S$. Then $f\circ\lambda$ is a closed maximal loxodromic curve on $\hat S$. Moreover, $\lambda$ contains no pole, so $f\circ\lambda$ contains no pole. Since every closed maximal loxodromic curve on a sphere-like surface is either a maximal parallel or a maximal meridian, and every maximal meridian contains both poles, $f\circ\lambda$ must be a maximal parallel. We have therefore proved the following proposition.

\begin{prop}\label{prop:map_to_parallel_and_meridian_for_Sphere}
Every loxodromic equivalence between sphere-like surfaces maps maximal meridians to maximal meridians and maximal parallels to maximal parallels.
\end{prop}

\begin{theorem}\label{thm:azimuth_consistent_sphere}
Let $S$ and $\hat S$ be sphere-like surfaces with sets of poles $P$
and $\hat P$, respectively. Set $S^*=S\setminus P$ and
$\hat S^*=\hat S\setminus\hat P$. A diffeomorphism
$f:S\to\hat S$ satisfying $f(P)=\hat P$ is azimuth-consistent if and
only if its differential on $S^*$ satisfies
\begin{align}
f_*\left(e_u\big|_p\right)
&=k(p)\hat e_u\big|_{f(p)},\label{eq:push_u_Sphere}\\
f_*\left(e_v\big|_p\right)
&=c k(p)\hat e_v\big|_{f(p)},\label{eq:push_v_Sphere}
\end{align}
where $p\in S^*$, $k$ is a smooth nonvanishing function on $S^*$,
and $c$ is a nonzero constant.
\end{theorem}

\begin{proof}
Suppose first that $f$ is azimuth-consistent. By
Proposition~\ref{prop:map_to_parallel_and_meridian_for_Sphere},
$f$ maps maximal meridians to maximal meridians and maximal parallels
to maximal parallels. Therefore, by the same argument as in
Theorem~\ref{thm:LP_and_OGmatrax}, its differential on $S^*$ has the
stated form.

Conversely, suppose that $f$ satisfies
\eqref{eq:push_u_Sphere} and \eqref{eq:push_v_Sphere}. Then both $f$
and $f^{-1}$ map loxodromes to loxodromes away from the poles, and
the directed azimuth of the image depends only on the directed
azimuth of the original loxodrome. The two equations also show that
both maps preserve meridional and parallel directions.

Since $f$ and $f^{-1}$ preserve meridional directions
and map poles to poles, maximal meridians are mapped
to maximal meridians. Maximal parallels are also mapped to maximal parallels. Every other maximal loxodromic curve has the two poles as its endpoints, and its image has the two poles of $\hat S$ as its endpoints. It follows that both $f$ and $f^{-1}$ map loxodromic curves to loxodromic curves.

The maximality argument in the converse part of
Theorem~\ref{thm:azimuth_consistent_paraboloid} shows that $f$
preserves maximal loxodromic curves. By
Proposition~\ref{prop:loxodromic-equivalence}, $f$ is a loxodromic
equivalence. Together with the preceding relation between the
directed azimuths, this shows that $f$ is azimuth-consistent.
\end{proof}

For an azimuth-consistent diffeomorphism between sphere-like
surfaces, the limiting argument used for paraboloid-like surfaces at
either pole gives $|c|=1$. Hence $c=\pm1$, and every
azimuth-consistent diffeomorphism between sphere-like surfaces is
conformal.

\subsection{Torus-Like Surfaces}\label{subsection:Torus-Like Surfaces}

The behavior of a loxodromic equivalence is considerably more complicated on torus-like surfaces. Specifically, these mappings might preserve neither parallels nor meridians. We provide an example of a loxodromic equivalence on a torus-like surface at the end of this subsection.

Let $S$ be a torus-like surface with closed profile curve
$\alpha=(r,0,z)$ defined on $[a,b]$, and let
$\Psi:\mathbb{R}^2\to S$ be its extended linearizing
parametrization, as given in
\eqref{eq:global_conformal_parametrization_on_torus}. Its period in the $w$-coordinate is $\Omega_w = \int_{a}^{b}\frac{1}{r(\tau)}\,d\tau$.

Consider a map $\mathcal{T}:\mathbb{R}^2\to\mathbb{R}^2$ defined by
\begin{equation}\label{eq:local_affine_map}
    \mathcal{T}\begin{pmatrix} u \\ w \end{pmatrix} = {D} M D^{-1} \begin{pmatrix} u \\ w \end{pmatrix} + \begin{pmatrix} b_1 \\ b_2 \end{pmatrix} ,
\end{equation}
where $D = \operatorname{diag}(2\pi, \Omega_w)$. Here, $b_1, b_2 \in \mathbb{R}$, and $M\in GL(2,\mathbb Z)$
is an integer matrix with determinant $\pm1$. Explicitly, we write $M = \begin{pmatrix}
m_{11} & m_{12} \\
m_{21} & m_{22}
\end{pmatrix} $ so that $D M D^{-1} = \begin{pmatrix} m_{11} & \frac{2\pi m_{12}}{\Omega_w} \\ \frac{{\Omega_w} m_{21}}{2\pi} &   m_{22} \end{pmatrix}.$

For any points $(u_1,w_1)$ and $(u_2,w_2)$ in $\mathbb{R}^2$ satisfying $\Psi(u_1,w_1) = \Psi(u_2,w_2)$, we have $(u_2,w_2) = (u_1 + 2\pi m,w_1 + \Omega_w n)$ for some integers $m, n$. It follows that
\begin{align*}
    \mathcal{T}(u_2,w_2)= \mathcal{T}(u_1,w_1)+\left(2\pi(m_{11}m+m_{12}n),\Omega_w(m_{21}m+m_{22}n)\right).
\end{align*}
Applying the periodicity of $\Psi$, we have
\begin{equation}\label{eq:forward_well_definedness}
    (\Psi\circ\mathcal{T})(u_1,w_1) = (\Psi\circ\mathcal{T})(u_2,w_2).
\end{equation}
In addition, since $M^{-1}\in GL(2,\mathbb Z)$, the map $D M^{-1}D^{-1}$ sends every vector in the period lattice $D\mathbb Z^2$ to another vector in $D\mathbb Z^2$. Therefore, the periodicity of $\Psi$ gives
\begin{equation}\label{eq:inverse_well_definedness}
    (\Psi\circ\mathcal{T}^{-1})(u_1,w_1) = (\Psi\circ\mathcal{T}^{-1})(u_2,w_2).
\end{equation}
Let $f:S \to S$ be a mapping satisfying
\begin{equation}\label{eq:f_definition}
    f(p) =  (\Psi\circ\mathcal{T})(u,w),
\end{equation}
where $(u,w)\in \mathbb{R}^2$ satisfies $p = \Psi(u,w)$. By \eqref{eq:forward_well_definedness}, $f$ is well-defined on $S$. We can define an inverse map $g:S \to S$ for $f$ by
$$
g(q) = (\Psi\circ\mathcal{T}^{-1})(\hat{u},\hat{w}),
$$
where $(\hat{u},\hat{w})\in \mathbb{R}^2$ satisfies $q = \Psi(\hat{u},\hat{w})$. Also, $g$ is well-defined by \eqref{eq:inverse_well_definedness}. Since
\begin{align*}
    (g\circ f)(p) &= (g\circ\Psi\circ \mathcal{T})(u,w) = (\Psi\circ\mathcal{T}^{-1}\circ\mathcal{T})(u,w) = p,\\
    (f\circ g)(p) &= (f\circ\Psi\circ \mathcal{T}^{-1})(u,w) = (\Psi\circ\mathcal{T}\circ\mathcal{T}^{-1})(u,w) = p,
\end{align*}
the mapping $f$ is a bijection, and its inverse is exactly $f^{-1}=g$.

\begin{prop}\label{prop:f_definition_is_diffeomorphism}
    A map $f:S \to {S}$ defined by \eqref{eq:f_definition} is a diffeomorphism.
\end{prop}
\begin{proof}
    Given $p \in S$, choose $(u,w) \in \mathbb{R}^2$ such that $p = \Psi(u,w)$. Let $q = f(p)$. Since $\Psi$ is a local diffeomorphism on $\mathbb{R}^2$, there exists a neighborhood $W \subset \mathbb{R}^2$ of $\mathcal{T}(u,w)$ such that the restriction $\Psi|_W: W \to \Psi(W)$ is a diffeomorphism. Since $\mathcal{T}$ is a diffeomorphism, $V = \mathcal{T}^{-1}(W)$ is a neighborhood of $(u,w)$. Similarly, since $\Psi$ is a local diffeomorphism, there exists a sufficiently small neighborhood $U \subseteq V$ of $(u,w)$ such that the restriction $\Psi|_U: U \to \Psi(U)$ is a diffeomorphism.

    Since $U \subseteq V$, we have $\mathcal{T}(U) \subseteq W$. The mapping $f$ can be expressed locally on the neighborhood $\Psi(U)$ as the composition $f|_{\Psi(U)}: \Psi(U) \to f(\Psi(U))$ given by
    $$
    f|_{\Psi(U)} = \Psi|_{\mathcal{T}(U)} \circ \mathcal{T}|_U \circ (\Psi|_U)^{-1}.
    $$
    Since $\Psi|_{\mathcal{T}(U)}$, $\mathcal{T}|_U$, and $(\Psi|_U)^{-1}$ are diffeomorphisms, the map $f|_{\Psi(U)}$ is also a diffeomorphism. Since $p$ is arbitrary, $f$ is a local diffeomorphism everywhere on $S$. Because $f$ is also bijective, we conclude that $f$ is a diffeomorphism.
\end{proof}

\begin{prop}\label{prop:loxo_diff_between_tori}
    The diffeomorphism $f:S\to S$ defined by \eqref{eq:f_definition} preserves maximal loxodromic curves.
\end{prop}
\begin{proof}
    Let $\lambda:\mathbb R\to S$ be a maximal loxodromic curve, and write $\lambda(t)=\Psi(c(t))$, where $c:\mathbb R\to\mathbb R^2$ has a straight line as its image. The curve $\mathcal T\circ c$ is also defined on $\mathbb R$ and has a straight line as its image. Hence $(f\circ\lambda)(t)=\Psi(\mathcal T(c(t)))$ is a loxodrome on $S$. Since $S$ is compact and $f$ is a diffeomorphism, the speed of $f\circ\lambda$ is bounded below by a positive constant. Therefore, an arc-length parametrization of $f\circ\lambda$ is defined on $\mathbb R$ and is maximal. Thus, $f$ preserves maximal loxodromic curves. By Proposition~\ref{prop:loxodromic-equivalence}, $f$ is a loxodromic equivalence.
\end{proof}
\begin{eg}\label{eg:torus_preservation_examples}
    Returning to Example~\ref{eg:linearizing_parameterization_of_torus}, let $f:S\to S$ be the map defined by \eqref{eq:f_definition} using a map of the form \eqref{eq:local_affine_map}. To illustrate that the diffeomorphism $f$ constructed in \eqref{eq:f_definition} preserves neither parallels nor meridians, let $\lambda_P(t) = \Psi(t,0)$ and $\lambda_M(t) = \Psi(0,t)$ be a parallel and a meridian, respectively, for $t \in \mathbb{R}$. Set
$$ M = \begin{pmatrix} 1 & 1 \\ 1 & 2 \end{pmatrix} \in GL(2, \mathbb{Z}) $$
    and $b = (0,0)^T$. Then we have $\mathcal{T}(t,0) = \left(t, \frac{\Omega_w}{2\pi}t\right)$, and $\mathcal{T}(0,t) = \left(\frac{2\pi}{\Omega_w}t, 2t\right)$. Since the $w$-coordinate of $(f\circ\lambda_P)(t) = \Psi\left(t, \frac{\Omega_w}{2\pi}t\right)$ is not constant, the image curve is not a parallel. Similarly, since the $u$-coordinate of $(f\circ\lambda_M)(t) = \Psi\left(\frac{2\pi}{\Omega_w}t, 2t\right)$ is not constant, the image curve is not a meridian. Thus, $f$ preserves neither parallels nor meridians.

In other cases, it is possible to choose a map $f$ that preserves both parallels and meridians. Set
$$ M = \begin{pmatrix} 1 & 0 \\ 0 & 1 \end{pmatrix} \in GL(2, \mathbb{Z}) $$
and $b = (1,2)^T$. Then we have $\mathcal{T}(t,0) = (t,0) + b = (t+1, 2)$, and $\mathcal{T}(0,t) = (0,t) + b = (1, t+2)$. The image curves are $(f\circ\lambda_P)(t) = \Psi(t+1, 2)$ and $(f\circ\lambda_M)(t) = \Psi(1, t+2)$. The first image has constant $w$-coordinate and is therefore
a parallel, while the second has constant $u$-coordinate
and is therefore a meridian. Thus, $f$ preserves both parallels and meridians.
\end{eg}

\bibliographystyle{IEEEtran}
\bibliography{refs}

\end{document}